\documentclass[english]{article}

\usepackage{geometry}
\usepackage{amsmath}
\usepackage{amsthm}
\usepackage{amssymb}
\usepackage{enumitem}
\usepackage{tabularx}
\usepackage{hyperref}
\usepackage{bm}
\usepackage{tabu}
\usepackage{blkarray}

\usepackage[noadjust]{cite}
\usepackage{hyperref}
\usepackage{mathtools}

\makeatletter

\newtheorem{theorem}{Theorem}[section]
\newtheorem{lemma}[theorem]{Lemma}
\newtheorem{proposition}[theorem]{Proposition}

\theoremstyle{definition}
\newtheorem{definition}[theorem]{Definition}
\newtheorem{example}[theorem]{Example}

\theoremstyle{remark}
\newtheorem{remark}[theorem]{Remark}

\numberwithin{equation}{section}

\usepackage{babel}

\newcommand{\bfH}{\mathbf{H}}
\newcommand{\bfh}{\mathbf{h}}
\newcommand{\bfI}{\mathbf{I}}
\newcommand{\bfA}{\mathbf{A}}
\newcommand{\bfB}{\mathbf{B}}
\newcommand{\bfM}{\mathbf{M}}
\newcommand{\bfS}{\mathbf{S}}
\newcommand{\bfv}{\mathbf{v}}

\newcommand{\bfx}{\mathbf{x}}
\newcommand{\bfy}{\mathbf{y}}

\newcommand{\bfg}{\mathbf{g}}

\newcommand{\rmH}{\mathrm{H}}
\newcommand{\rmA}{\mathrm{A}}

\newcommand{\rmT}{\mathrm{T}}

\newcommand{\calF}{{\mathcal{F}}_{\mu,L}}

\DeclareMathOperator{\spann}{span}

\begin{document}

\title{A Proof of the Exact Convergence Rate of Gradient Descent\footnote{Part~I (\texttt{arXiv:2412.04435v1}) and Part~II (\texttt{arXiv:2412.04427v1}) of this paper are merged into this version.}}

\author{Jungbin Kim}
\date{}

\maketitle

\begin{abstract}
We prove the exact worst-case convergence rate of gradient descent for smooth strongly convex optimization on $\mathbb{R}^d$. Concretely, assuming that the objective function $f$ is $\mu$-strongly convex and $L$-smooth, we identify the smallest possible value of $\tau$ for which the inequality $f(x_{N})-f_{*}\leq\tau\|x_{0}-x_{*}\|^{2}$ always holds. The result was previously conjectured by Drori and Teboulle for the case $\mu=0$, and by Taylor, Hendrickx, and Glineur for the case $\mu>0$.
\end{abstract}

\section{Introduction}
\label{sec:intro}

We study the problem of minimizing a differentiable function $f$ on $\mathbb{R}^d$, that is,
\begin{equation}
	\label{eq:opt_prob}
	\min_{x\in\mathbb{R}^{d}}\; f(x).
\end{equation}
The algorithm we consider is gradient descent (GD) with a fixed stepsize $\gamma>0$. Given an initial point $x_0 \in \mathbb{R}^d$, the algorithm updates the iterates as
\begin{equation}
	\tag{\textup{GD}}
	\label{eq:gd}
	x_{k+1}=x_k - \gamma\nabla f(x_k)\quad\textup{for }k=0,1,\ldots.
\end{equation}
In this paper, we answer the following open question in convex optimization: ``How quickly does gradient descent converge to the minimum?'' This informal question will be mathematically formulated in the next two subsections.

\subsection{Assumptions}
\label{sec:assumptions}

In this subsection, we summarize the standard assumptions which we will follow throughout the paper. We assume that $f$ has a minimizer $x_*$, which may not be unique. When dealing with quantities depending on $x_*$, we consider $x_*$ to be chosen arbitrarily among them. We denote $f_*=f(x_*)$. 

It is common to assume that $f$ is convex in order to make the problem \eqref{eq:opt_prob} tractable. To ensure that \ref{eq:gd} with a predetermined fixed stepsize always converges, it is usually assumed that the gradient of $f$ is Lipschitz continuous with a known Lipschitz constant. For our purposes, the following working definitions suffice to cover these two assumptions: Given parameters $-\infty<\mu<L<\infty$, we assume that
\begin{itemize}
	\item $f$ is $\mu$-strongly convex, that is, for all $p, q \in \mathbb{R}^d$, we have
	\begin{equation}
		\label{eq:strong-convexity-def}
		0 \geq f(p) - f(q) + \langle \nabla f(p), q - p \rangle + \frac{\mu}{2} \| p - q \|^2;
	\end{equation}
	\item $f$ is $L$-smooth, that is, for all $p, q \in \mathbb{R}^d$, we have
	\begin{equation}
		\label{eq:smoothness-def}
		0 \leq f(p) - f(q) + \langle \nabla f(p), q - p \rangle + \frac{L}{2} \| p - q \|^2.
	\end{equation}
\end{itemize}
Throughout the paper, we always assume $L>0$. Depending on the context, we sometimes require $\mu \geq 0$. Keep in mind that we allow $\mu<0$ unless specified otherwise. When $\mu=0$, the notion of $\mu$-strong convexity is equivalent to the standard notion of convexity. When $f$ is convex, it is known that the $L$-smoothness of $f$ is equivalent to the $L$-Lipschitz continuity of $\nabla f$ (see, for example, \cite[Thm.~2.1.5]{nesterov2018lectures}). We denote by $\calF$ the space of $\mu$-strongly convex and $L$-smooth functions on $\mathbb{R}^d$ which have a minimizer. The ratio $\mu/L$ is denoted by $\kappa$. The following property is well-known (see, for example, \cite[Thm.~2.1.5]{nesterov2018lectures}, \cite[Thm.~A.3]{daspremont2021acceleration}, and \cite[Thm.~3.2]{rotaru2024exact}).
\begin{proposition}
	\label{prop:int-ineq}
	When $f\in\calF$, the \emph{interpolation inequality}
	\begin{equation}
		\label{eq:interpolation_ineqq}
		\begin{aligned}
			0 & \geq f(p)-f(q)+\left\langle \nabla f(p),q-p\right\rangle \\
			& \quad + \frac{1}{2L}\left\Vert \nabla f(p)-\nabla f(q)\right\Vert ^{2}+\frac{\mu L}{2(L-\mu)}\left\Vert p-q-\frac{1}{L}\left(\nabla f(p)-\nabla f(q)\right)\right\Vert ^{2}
		\end{aligned}
	\end{equation}
	holds for all $p,q\in\mathbb{R}^d$.
\end{proposition}

\subsection{History of the problem}
\label{sec:history}

The study of the convergence properties of gradient descent has a long history. The method of gradient descent can be traced back to \cite{cauchy1847methode,hadamard1908memoire}. An early work on the convergence of gradient descent is \cite{curry1944method}, where neither the convexity of $f$ nor the Lipschitz continuity of $\nabla f$ is assumed. As a result, line search (i.e., finding an appropriate stepsize $\gamma>0$ at each iteration) was necessary, and the rate of convergence was not shown. In the following decades, the convergence properties of \ref{eq:gd} became well-studied as the standard assumptions in Section~\ref{sec:assumptions} were adopted. In particular, it was shown in \cite[Thm.~6]{polyak1963gradient} that when $f\in\calF$ with $\mu\geq 0$, the iterates of \ref{eq:gd} with $\gamma\in(0,2/L)$ satisfy the inequality
\begin{equation}
	\label{eq:polyak}
	\left\Vert x_{N}-x_{*}\right\Vert \leq\max\left\{ |1-\gamma\mu|,|1-\gamma L|\right\} ^{N}\left\Vert x_{0}-x_{*}\right\Vert .
\end{equation}
When $\mu > 0$, this inequality not only shows that $x_N\to x_*$ as $N \to \infty$, but also provides an upper bound on the rate at which the distance $\|x_N - x_*\|$ decreases. We refer to the ratio which we aim to bound above ($\|x_{N}-x_{*}\|/\|x_{0}-x_{*}\|$ in this case) as the \emph{performance criterion}.

From a practical perspective, it is more reasonable to measure the quality of $x_N$ as an approximate solution for the optimization problem \eqref{eq:opt_prob} using $f(x_N)-f_*$ rather than $\|x_N-x_*\|$.\footnote{Another reason to consider the performance criterion $(f(x_N)-f_*)/\|x_0-x_*\|^2$ is that when $\mu=0$, one cannot prove meaningful bounds for the performance criteria $\|x_N-x_*\|/\|x_0-x_*\|$ and $(f(x_N)-f_*)/(f(x_0)-f_*)$. In those cases, $\|x_N-x_*\|\leq\|x_0-x_*\|$ and $f(x_N)-f_*\leq f(x_0)-f_*$ are the best possible bounds (see \cite{taylor2018exact}), indicating that \ref{eq:gd} yields no improvement in these respects.} For example, the following bounds can be found in \cite[\S 2.1.5]{nesterov2018lectures} (see also \cite[\S3]{bubeck2015convex} and \cite[\S4]{daspremont2021acceleration}).
\begin{proposition}
	\label{prop:poor-rate}
	When $f\in \mathcal{F}_{0,L}$, the iterates of \ref{eq:gd} with $\gamma=1/L$ satisfy
	\begin{equation}
		\label{eq:gd-poor-rate}
		f(x_{N})-f_*\leq\frac{2L}{N+4}\left\Vert x_{0}-x_{*}\right\Vert ^{2}.
	\end{equation}
	When $f\in \mathcal{F}_{\mu,L}$ with $\mu\geq 0$, the iterates of \ref{eq:gd} with $\gamma=2/(\mu+L)$ satisfy
	\begin{equation}
		\label{eq:gd-poor-rate2}
		f(x_{N})-f_*\leq\frac{L}{2}\left(\frac{L-\mu}{L+\mu}\right)^{2N}\left\Vert x_{0}-x_{*}\right\Vert ^{2}.
	\end{equation}
\end{proposition}
Here, the performance criterion is set as $(f(x_N)-f_*)/\|x_0-x_*\|^2$, which is perhaps the most widely studied performance criterion in the convex optimization literature. A natural question is whether the inequalities \eqref{eq:gd-poor-rate} and \eqref{eq:gd-poor-rate2} can be tightened. This leads to the problem of identifying the value
\begin{equation}
	\label{eq:exact-rate-def}
	\sup\left\{ \frac{f(x_{N})-f_*}{\left\Vert x_{0}-x_{*}\right\Vert ^{2}}:f\in\calF,\ x_0\in\mathbb{R}^d,\ x_{k}\textup{ generated by \ref{eq:gd}}\right\} 
\end{equation}
for given $\gamma$, $N$, $\mu$, and $L$. Since $f(x_k)$ may not converge when $\gamma\notin(0,2/L)$,\footnote{Consider the objective function $f(x)=(L/2)\|x\|^2$, which is in $\calF$.} we assume $\gamma\in(0,2/L)$. Based on numerical evidence, a closed-form expression of the value \eqref{eq:exact-rate-def} was conjectured in \cite{drori2014performance} for the case $\mu=0$. Later, this conjecture was generalized in \cite{taylor2017smooth} to encompass the case $\mu > 0$. The conjectured result is shown in the following theorem, which we will prove.
\begin{theorem}[{\cite[Conj.~1]{drori2014performance}, \cite[Conj.~2]{taylor2017smooth}}]
	\label{conj:gd}
	When $f\in \mathcal{F}_{\mu,L}$ with $\mu\geq 0$, the iterates of \ref{eq:gd} with $\gamma\in(0,2/L)$ satisfy
	\begin{equation}
		\label{eq:gd-rate}
		f(x_{N})-f_{*}\leq\max\left\{ \frac{1}{1+\gamma L\sum_{j=1}^{2N}(1-\gamma\mu)^{-j}},(1-\gamma L)^{2N}\right\} \frac{L\left\Vert x_{0}-x_{*}\right\Vert ^{2}}{2}.
	\end{equation}
\end{theorem}
Note that Theorem~\ref{conj:gd} only provides an upper bound of \eqref{eq:exact-rate-def}. To identify the value \eqref{eq:exact-rate-def}, we also need the matching lower bound. This is done by the following proposition.
\begin{proposition}[{\cite{drori2014performance}, \cite[\S4.1.2]{taylor2017smooth}}]
	\label{prop:lower-f}
	There is a pair $(f,x_0)\in\calF\times\mathbb{R}^d$ for which \eqref{eq:gd-rate} holds with equality.
\end{proposition}
To the best of our knowledge, the upper bound (Theorem~\ref{conj:gd}) has only been partially proved for the case where $\mu=0$ and $\gamma \in (0, 1/L]$ in \cite[Thm.~2]{drori2014performance} and \cite[Thm.~2]{teboulle2023elementary}, while the lower bound (Proposition~\ref{prop:lower-f}) was completely proved in \cite{drori2014performance} for the case $\mu=0$,\footnote{See Theorem~3.2 of the preprint version (\texttt{arXiv:1206.3209v1}) of their paper.} and in \cite{taylor2017smooth} for the case $\mu>0$.

Although the connection to Theorem~\ref{conj:gd} was not clear, a recent work \cite{rotaru2024exact} made notable progress by proving the exact convergence rate of \ref{eq:gd} for another performance criterion $\Vert \nabla f(x_N)\Vert^2 / (f(x_0)-f_*)$. As in the previous case, both the upper and lower bounds are needed. Their results are shown below.
\begin{theorem}[{\cite[Thm.~2.2]{rotaru2024exact}}]
	\label{thm:gd}
	When $f\in \mathcal{F}_{\mu,L}$ with $\mu\geq 0$, the iterates of \ref{eq:gd} with $\gamma\in(0,2/L)$ satisfy
	\begin{equation}
		\label{eq:gd-rate2}
		\frac{\left\Vert \nabla f(x_{N})\right\Vert ^{2}}{2L}\leq\max\left\{ \frac{1}{1+\gamma L\sum_{j=1}^{2N}(1-\gamma\mu)^{-j}},(1-\gamma L)^{2N}\right\} \left(f(x_{0})-f_{*}\right).
	\end{equation}
\end{theorem}
\begin{proposition}[{\cite[Prop.~5.2]{rotaru2024exact}}]
	There is a pair $(f,x_0)\in\calF\times\mathbb{R}^d$ for which \eqref{eq:gd-rate2} holds with equality.
\end{proposition}
In another recent work \cite{grimmer2024strengthened}, it was conjectured that Theorem~\ref{conj:gd} for the optimal stepsize $\gamma=\gamma^*(N,\mu,L)$ (see Section~\ref{sec:proof_special}) can be proved using the PEP methodology (see Section~\ref{sec:review_drori}) with a structured set of multipliers $\{\lambda_{i,j}\}_{i,j\in I_N^*}$. Although we prove Theorem~\ref{thm:F2} in a different way from the conjectured one, our proof and their conjectured proof share similarities. See Remark~\ref{rmk:grimmer} for details.

\subsection{Previously known partial proof}
\label{sec:review_drori}

In this subsection, we review the partial proof of Theorem~\ref{conj:gd} presented in \cite{drori2014performance}, which is based on the methodology of \emph{performance estimation problem} (PEP).\footnote{By definition, a performance estimation problem (PEP) is an optimization problem of finding the value \eqref{eq:exact-rate-def} (or its variants). The analysis in this subsection was originally derived in \cite{drori2014performance} by taking the Lagrangian dual of a certain relaxation of PEP. We will refer to any convergence analysis involving the verification of inequalities like \eqref{eq:target} as a {PEP analysis}.} It is assumed that $\mu=0$.

Introduce the index set $I_N^{*}:=\{0,\ldots,N,*\}$. Suppose that there is a set of nonnegative \emph{multipliers} $\{\lambda_{i,j}\}_{i,j\in I_{N}^{*}}$ such that the inequality
\begin{equation}
	\label{eq:target}
	\begin{aligned}
		& f_{N}-f_{*}-\max\left\{ \frac{1}{1+2N\gamma L},(1-\gamma L)^{2N}\right\} \frac{L\left\Vert x_{0}-x_{*}\right\Vert ^{2}}{2}\\
		& \qquad\leq\sum_{i,j\in I_{N}^{*}}\lambda_{i,j}\left(f_{j}-f_{i}+\left\langle g_{j},x_{i}-x_{j}\right\rangle +\frac{1}{2L}\left\Vert g_{i}-g_{j}\right\Vert ^{2}\right)
	\end{aligned}
\end{equation}
holds for all [families $\{x_i\}_{i\in I_N^*}$ in $\mathbb{R}^d$, $\{g_i\}_{i\in I_N^*}$ in $\mathbb{R}^d$, and $\{f_i\}_{i\in I_N^*}$ in $\mathbb{R}$ satisfying $g_*=0$ and $x_{k+1}=x_k-\gamma g_k$ for $k=0,\ldots,N-1$]. Then, we can put $x_i\leftarrow x_i$,\footnote{\label{footnote:abuse}We abused notation here. The left-hand $x_k$ refer to arbitrary points in the previous sentence, while the right-hand $x_k$ denote the iterates of \ref{eq:gd} and the given minimizer $x_*$.} $g_i\leftarrow\nabla f(x_i)$, and $f_i\leftarrow f(x_i)$ for all $i\in I_N^*$, and then the right-hand side of \eqref{eq:target} is nonpositive by Proposition~\ref{prop:int-ineq}. Thus, the inequality \eqref{eq:target} implies the desired inequality \eqref{eq:gd-rate}. 

Therefore, to prove Theorem~\ref{conj:gd}, it suffices to set appropriate multipliers $\lambda_{i,j}$ and then verify that the inequality \eqref{eq:target} holds under the stated conditions. An advantage of this approach is that the task of finding working multipliers $\lambda_{i,j}$ can be formulated as a semidefinite program, and thus one can use numerical solvers to confirm their existence and even infer their closed-form expressions. A difficulty that arises here is that feasible solutions $\{\lambda_{i,j}\}_{i,j\in I_N^*}$ to this semidefinite program are not unique. In other words, there is some freedom in choosing $\{\lambda_{i,j}\}_{i,j\in I_N^*}$. This causes numerical solvers to yield inconsistent solutions. To resolve this issue (although they did not explicitly mention the reason), the authors of \cite{drori2014performance} imposed an additional constraint 
\begin{equation}
	\label{eq:relaxed-f}
	\lambda_{i,j}=0\quad\textup{for all }(i,j)\notin\left\{ (k,k+1),(*,k)\right\} _{0\leq k\leq N-1}\cup\left\{ (*,N)\right\} .
\end{equation}
This constraint removes the freedom in selecting $\{\lambda_{i,j}\}_{i,j\in I_N^*}$, enabling them to infer a closed-form expression of the solution $\{\lambda_{i,j}\}_{i,j\in I_N^*}$. This strategy was successful for the case $\gamma\in(0,1/L]$, leading to a partial proof of Theorem~\ref{conj:gd}. However, they were not able to extend their proof to the case $\gamma\in(1/L,2/L)$, since feasible solutions $\{\lambda_{i,j}\}_{i,j\in I_N^*}$ satisfying \eqref{eq:relaxed-f} do not always exist in this case (this is suggested by the numerical evidence in \cite[Table~1]{taylor2017smooth}).

\subsection{Our strategy}
\label{sec:our_strategy}

In this paper, we prove Theorem~\ref{conj:gd} for any strong convexity parameter $\mu\in[0,L)$ and any stepsize $\gamma\in(0,2/L)$. We use the PEP methodology in the previous subsection, but even for the case where $\mu=0$ and $\gamma\in(0,1/L]$, we choose a different set of multipliers $\{\lambda_{i,j}\}_{i,j\in I_N^*}$ from the one used in \cite{drori2014performance}. 

A crucial ingredient of our proof is observing the phenomenon that some convergence analyses for a first-order optimization algorithm $\bfH$ with respect to the performance criterion $(f(x_N)-f_*)/\Vert x_0-x_*\Vert^2$ have \emph{dual} convergence analyses for the \emph{anti-transposed} algorithm $\bfH^{\rmA}$ with respect to the performance criterion $\Vert \nabla f(x_N)\Vert^2 / (f(x_0)-f_*)$. See \cite[Thm.~3]{kim2024convergence} and \cite[Thms.~1,2]{kim2024time} for a recent theory on this phenomenon. To compare our approach for obtaining a proof of Theorem~\ref{conj:gd} with the one in \cite{drori2014performance}, we handle the nonuniqueness of feasible solutions $\{\lambda_{i,j}\}_{i,j\in I_N^*}$ by removing the freedom in their choice within the {dual} analysis.\footnote{We empirically observed that our proof strategy fully removes the freedom in choosing the multipliers. We will not prove this because it is not essential to our main objective.} 

The structure of our paper is outlined as follows. In Section~\ref{sec:4}, we develop a theory for the aforementioned phenomenon. Our theory encompasses the existing theory and is supported by a new proof. In this correspondence, proving Theorem~\ref{thm:gd} will mirror the task of proving Theorem~\ref{conj:gd}. In Section~\ref{sec:pep}, we provide a proof of Theorem~\ref{thm:gd}, which differs from the previous one presented in \cite{rotaru2024exact}. In Section~\ref{sec:gd}, we combine the two results to complete the proof of Theorem~\ref{conj:gd}.

Our proof involves tedious symbolic computations. Rather than presenting written calculations, we provide MATLAB codes for verifying the computations. The codes are available at the following GitHub repository:
\begin{center}
	\url{https://github.com/jungbinkim1/GD-Exact-Rate}.
\end{center}

\section{Preliminaries and notation}
\label{sec:notation}

We mostly follow the notations in \cite{rotaru2024exact}. We define $\rho=1-\gamma L$ and $\eta=1-\gamma \mu $. For $k=1,2,\ldots$, we define functions $E_k(\cdot)$, $F_k(\cdot)$, and $T_k(\cdot,\cdot)$ as
\begin{align*}
	E_{k}(x) & {\textstyle =\begin{cases}
			\sum_{j=1}^{2k}x^{-j} & \textup{if }x\neq0,\\
			\infty & \textup{if }x=0,
	\end{cases}}\\
	F_{k}(x) & {\textstyle =\sum_{j=1}^{k}x^{j},}\\
	T_{k}(\rho,\eta) & {\textstyle =E_{k}(\eta)-E_{k}(\rho).}
\end{align*}
See Table~\ref{table:1} for a rough description of the function $E_k$. We introduce the index sets $I_N=\{0,\ldots,N\}$ and $I_N^*=I_N\cup\{*\}$. For convenience, we set all multipliers $\lambda_{i,j}$ (or $\nu_{i,j}$) with $i=j$ to zero throughout the paper.

We denote the space of $N$-tuples of points in $\mathbb{R}^d$ by $(\mathbb{R}^d)^N$. We represent elements of this space as column vectors, that is, $[p_{1}, \ldots, p_{N}]^{\rmT}$ with $p_k \in \mathbb{R}^d$ for $k = 1, \ldots, N$. To point out trivial facts, we note that the space $(\mathbb{R}^d)^N$ is equipped with the inner product $\langle[p_{1},\ldots,p_{N}]^{\rmT},[q_{1},\ldots,q_{N}]^{\rmT}\rangle=\sum_{k=1}^{N}\langle p_{k},q_{k}\rangle$, and that $N \times N$ matrices can be viewed as operators on $(\mathbb{R}^d)^N$.

\subsection{Properties of interpolation inequality}

The following proposition shows that the interpolation inequality \eqref{eq:interpolation_ineqq} tightens as $\mu$ increases or $L$ decreases.
\begin{proposition}
	\label{prop:ineq-stronger}
	Fix $p$, $q$, $\nabla f(p)$, $\nabla f(q)$, $f(p)$, and $f(q)$. Denote the right-hand side of \eqref{eq:interpolation_ineqq} as $Q(\mu,L)$. Then, we have
	\begin{enumerate}[label=\textnormal{(\roman*)},ref=\roman*]
		\item \label{item:mu-compared} $Q(\mu,L)\geq Q(\mu',L)$ whenever $\mu>\mu'$;
		\item \label{item:L-compared} $Q(\mu,L)\geq Q(\mu,L')$ whenever $L<L'$.
	\end{enumerate}
\end{proposition}
\begin{proof}
	These follow from
	\begin{align*}
		\frac{\partial Q}{\partial\mu} & =\frac{L^{2}}{2(L-\mu)^{2}}\left\Vert p-q-\frac{1}{L}\left(\nabla f(p)-\nabla f(q)\right)\right\Vert ^{2}\geq 0,\\
		\frac{\partial Q}{\partial L} & =-\frac{\mu^{2}}{2(L-\mu)^{2}}\left\Vert \left(p-q\right)-\frac{1}{\mu}\left(\nabla f(p)-\nabla f(q)\right)\right\Vert ^{2}\leq 0.
	\end{align*}
\end{proof}
Since $\nabla f(x_*)=0$, the inequality \eqref{eq:interpolation_ineqq} with $p=x_*$ can be written as
\[
0\geq f_*-f(q)+\frac{1}{2L}\left\Vert \nabla f(q)\right\Vert ^{2}+\frac{\mu L}{2(L-\mu)}\left\Vert x_{*}-q+\frac{1}{L}\nabla f(q)\right\Vert ^{2}.
\]
Thus, when $\mu\geq0$, we have $f(q)-\frac{1}{2L}\|\nabla f(q)\|^2\geq f_*$. The following proposition shows that this inequality remains valid even when $\mu<0$.
\begin{proposition}
	\label{prop:trivial}
	If $f$ is $L$-smooth, then the inequality $f(x)-\frac{1}{2L}\|\nabla f(x)\|^2\geq f_*$ holds for all $x\in\mathbb{R}^d$.
\end{proposition}
\begin{proof}
	By putting $p=x$ and $q=x-\frac{1}{L}\nabla f(x)$ in the inequality \eqref{eq:smoothness-def}, we obtain $f(x-\frac{1}{L}\nabla f(x))\leq f(x)-\frac{1}{2L}\|\nabla f(x)\|^2$. Combine this with $f_*\leq f(x-\frac{1}{L}\nabla f(x))$.
\end{proof}

\begin{table}
	\caption{The behavior of the the function $E_k$ for any fixed $k$.}
	\label{table:1}
	\centering
	\begin{tabular}{|c||c|c|c|c|c|} 
		\hline
		$x$ & $-1$ & $\cdots$ & $0$ & $\cdots$ & $\infty$ \\
		\hline
		$E_k(x)$ & $0$  & $\nearrow$ & $\infty$ & $\searrow$ & $0$ \\
		\hline
	\end{tabular}
\end{table}

\subsection{Optimiation algorithms as lower-triangular matrices}
\label{sec:h-matrix-rep}

Consider first-order methods where each iterate $x_{k}$ is in $x_0 + \spann\{\nabla f(x_0),\ldots,\nabla f(x_{k-1})\}$. For such a method, there is an $N\times N$ lower-triangular matrix $\bfH$ such that
\begin{equation}
	\label{eq:fsfom}
	\left[\begin{array}{c}
		x_{1}-x_{0}\\
		\vdots\\
		x_{N}-x_{N-1}
	\end{array}\right]=-\frac{1}{L}\left[\begin{array}{ccc}
		\bfH_{1,1}\\
		\vdots & \ddots\\
		\bfH_{N,1} & \cdots & \bfH_{N,N}
	\end{array}\right]\left[\begin{array}{c}
		g_{0}\\
		\vdots\\
		g_{N-1}
	\end{array}\right],
\end{equation}
where $g_k=\nabla f(x_k)$ for $k=0,\ldots,N-1$. In particular, \ref{eq:gd} can be expressed as \eqref{eq:fsfom} with $\bfH = (\gamma L) \bfI_N$, where $\bfI_N$ is the $N\times N$ identity matrix. Conversely, any lower-triangular matrix $\bfH$ gives an implementable algorithm \eqref{eq:fsfom}. Therefore, we can identify the algorithm \eqref{eq:fsfom} with the lower-triangular matrix $\bfH$. 

We can reparametrize \eqref{eq:fsfom} so that the iterates on the left-hand side become $x_k^+=x_k-\frac{1}{L}g_k$ instead of $x_k$. Define an $(N+1)\times (N+1)$ matrix $\tilde{\bfH}$ as
\begin{equation}
	\label{eq:h-tilde}
	\tilde{\bfH}=\left[\begin{array}{ccccc}
		1\\
		\bfH_{1,1}-1 & 1\\
		\vdots & \vdots & \ddots\\
		\bfH_{N-1,1} & \bfH_{N-1,2} & \cdots & 1\\
		\bfH_{N,1} & \bfH_{N,2} & \cdots & \bfH_{N,N}-1 & 1
	\end{array}\right].
\end{equation}
Then, \eqref{eq:fsfom} can be equivalently written as
\begin{equation}
	\label{eq:fsfom-tilde}
	\bfx=-\frac{1}{L}\tilde{\bfH}\bfg,
\end{equation}
where $\bfx=[x_{0}^{+}-x_{0},x_{1}^{+}-x_{0}^{+},\ldots,x_{N}^{+}-x_{N-1}^{+}]^{\rmT}$ and $\bfg=[g_{0},\ldots,g_{N}]^{\rmT}$.

\subsection{Anti-transpose of optimization algorithm}
\label{sec:at}

For an $N\times N$ matrix $\bfH$, its \emph{anti-diagonal transpose} (or shortly \emph{anti-transpose}) $\bfH^{\rmA}$ is obtained by reflecting $\bfH$ over its anti-diagonal, that is, $\bfH^{\rmA}_{i,j}=\bfH_{N+1-j,N+1-i}$ for $i,j=1,\ldots,N$. When $\bfH$ is lower-triangular, the first-order algorithm associated with $\bfH^{\rmA}$ can be written as
\begin{equation}
	\label{eq:fsfom-dual}
	\left[\begin{array}{c}
		y_{1}-y_{0}\\
		\vdots\\
		y_{N}-y_{N-1}
	\end{array}\right]=-\frac{1}{L}\left[\begin{array}{ccc}
		\bfH_{N,N}\\
		\vdots & \ddots\\
		\bfH_{N,1} & \cdots & \bfH_{1,1}
	\end{array}\right]\left[\begin{array}{c}
		g_0\\
		\vdots\\
		g_{N-1}
	\end{array}\right],
\end{equation}
where $g_k=\nabla f(y_k)$ for $k=0,\ldots,N-1$. We denoted the iterates by $y_k$ instead of $x_k$ to prevent future confusion. In particular, the anti-transpose of the matrix $(\gamma L) \bfI_N$ is itself. Perhaps, the observation of the anti-transpose relationship between first-order optimization algorithms dates back to \cite[\S6.6]{kim2021optimizing}. Note that $\widetilde{\bfH^{\rmA}} = \tilde{\bfH}^{\rmA}$ holds, where the tilde in the notation is defined as in \eqref{eq:h-tilde}. Thus, with $y_k^+=y_k-\frac{1}{L}g_k$, \eqref{eq:fsfom-dual} can be equivalently written as
\begin{equation}
	\label{eq:fsfom-dual-tilde}
	\bfy=-\frac{1}{L}\tilde{\bfH}^{\rmA}\bfg,
\end{equation}
where $\bfy=[y_{0}^{+}-y_{0},y_{1}^{+}-y_{0}^{+},\ldots,y_{N}^{+}-y_{N-1}^{+}]^{\rmT}$ and $\bfg=[g_{0},\ldots,g_{N}]^{\rmT}$.

The anti-transpose can be defined using the language of operators as follows.
\begin{definition}
	\label{def:at}
	Let $\rmH$ be an operator on $(\mathbb{R}^d)^N$. Then, its anti-transpose is the operator $\rmH^{\rmA}$ on $(\mathbb{R}^d)^N$ that satisfies $$\langle \rmH x, \rmT y \rangle = \langle \rmT x, {\rmH}^{\rmA} y \rangle$$ for all $x,y\in(\mathbb{R}^d)^N$, where $\rmT$ is the time reversal operator on $(\mathbb{R}^d)^N$ defined as 
	\[
	\rmT\left([x_{1},\ldots,x_{N}]^{\rmT}\right)=[x_{N},\ldots,x_{1}]^{\rmT}.
	\]
\end{definition}
Then, by definition, we have $\rmH^{\rmA} = \rmT \rmH^{\rmT} \rmT$ (this identity can be found in \cite[\S1]{golyshev2007fuchsian}). When $\rmH$ is invertible, from the identity $(\rmH^{\rmT})^{-1}=(\rmH^{-1})^{\rmT}$, it follows that 
\[
\left(\rmH^{\rmA}\right)^{-1}=\left(\rmT\rmH^{\rmT}\rmT\right)^{-1}=\left(\rmT\right)^{-1}\left(\rmH^{\rmT}\right)^{-1}\left(\rmT\right)^{-1}=\rmT\left(\rmH^{-1}\right)^{\rmT}\rmT=\left(\rmH^{-1}\right)^{\rmA}.
\]
This justifies denoting this operator as ${\rmH}^{-\rmA}$.

\section{Equivalence between convergence analyses for minimizing $(f(x_N)-f_*)/\|x_{0}-x_{*}\|^{2}$ and minimizing $\|\nabla f(x_{N})\|^{2}/(f(x_0)-f_*)$}
\label{sec:4}

In this section, we show that there is an equivalence between some PEP analyses for the algorithm $\bfH$ with the performance criterion $(f(x_N)-f_*) / \|x_{0}-x_{*}\|^{2}$ and some PEP analyses for the algorithm $\bfH^{\rmA}$ with the performance criterion ${\|\nabla f(x_{N})\|^{2}} / {(f(x_0)-f_*)}$.
\subsection{PEP analysis for performance criterion ${(f(x_N)-f_*)} / {\|x_{0}-x_{*}\|^{2}}$}
\label{sec:F2}

In this subsection, we present a PEP analysis aimed at proving convergence guarantees in the form $f(x_N)-f_*\leq\frac{\tau}{2}\|x_0-x_*\|^2$. Our analysis differs from the one in \cite{drori2014performance} by using $\{x_{k+1}^+ - x_k^+\}$ instead of $\{g_k\}$ as the basis for the quadratic form,\footnote{\label{footnote:see-psd}This will become clearer in Section~\ref{sec:detail_pep}, where we present a PEP analysis in terms of the positive semidefiniteness of a certain matrix.} although both analyses aim at the same performance criterion.
\begin{theorem}
	\label{thm:F2}
	Given a set of nonnegative multipliers $\{\lambda_{i,j}\}_{i,j\in I_N^*}$, the inequality
	\begin{equation}
		\label{eq:F2-rate}
		\begin{aligned}
			& f_{N}-f_{*}-\frac{\tau}{2}\left\Vert x_{0}-x_{*}\right\Vert ^{2}\\
			& \qquad\leq\sum_{i,j\in I_{N}^{*}}\lambda_{i,j}\bigg[f_{j}-f_{i}+\left\langle g_{j},x_{i}-x_{j}\right\rangle +\frac{1}{2L}\left\Vert g_{i}-g_{j}\right\Vert ^{2}\\
			& \qquad\qquad\qquad\qquad\quad+\frac{\mu L}{2(L-\mu)}\left\Vert x_{i}-x_{j}-\frac{1}{L}\left(g_{i}-g_{j}\right)\right\Vert ^{2}\bigg]
		\end{aligned}
	\end{equation}
	holds for all [families $\{x_i\}_{i\in I_N^*}$ in $\mathbb{R}^d$, $\{g_i\}_{i\in I_N^*}$ in $\mathbb{R}^d$, and $\{f_i\}_{i\in I_N^*}$ in $\mathbb{R}$ satisfying $g_*=0$ and \eqref{eq:fsfom}] if and only if the following two conditions are satisfied:
	\begin{enumerate}[label=\textnormal{(F\arabic*)},ref=F\arabic*]
		\item \label{item:F2-1} \[
		\sum_{i\in I_N^*}\lambda_{i,k}-\sum_{j\in I_N^*}\lambda_{k,j}=\begin{cases}
			1 & \textup{if }k=N,\\
			-1 & \textup{if }k=*,\\
			0 & \textup{otherwise}.
		\end{cases}
		\]
		\item \label{item:F2-2} For any points $x_{0},x_{0}^{+},\ldots,x_{N}^{+},x_*^+,y_{0},y_{0}^{+},\ldots,y_{N}^{+},y_*^+\in\mathbb{R}^d$ satisfying
		\begin{equation}
			\label{eq:xmapy2}
			\begin{aligned}
				y_{k}^{+}-y_{k-1}^{+}& =\sum_{i\in I_N^*}\lambda_{i,N-k}\left(x_{i}^{+}-x_{N-k}^{+}\right)\quad\textup{for }k=1,\ldots,N,\\
				\frac{y_{0}^{+}+y_{*}^+}{2}-y_{N}^{+}& =\frac{x_{N}^{+}-x_{*}^+}{2}+\sum_{i\in I_N^*}\lambda_{i,*}\left(x_{i}^{+}-x_{*}^+\right),
			\end{aligned}
		\end{equation}
		the quantity
		\begin{equation}
			\label{eq:F-PEP-functional2}
			\begin{aligned}
				S & =-L\left\langle \tilde{\bfH}^{-1}\bfx,\rmT\bfy\right\rangle +\frac{L}{2}\left\Vert y_{0}-\frac{y_{0}^{+}+y_{*}^+}{2}-\frac{x_{N}^{+}-x_{*}^+}{2}\right\Vert ^{2}\\
				& \quad-\frac{\mu L}{L-\mu}\bigg[\sum_{j=1}^{N}\left\langle y_{k}^{+}-y_{k-1}^{+},x_{N-k}^{+}\right\rangle  \\
				& \qquad\qquad\qquad +\left\langle \frac{y_{0}^{+}+y_{*}^+}{2}-y_{N}^{+},x_{*}^+\right\rangle  +\left\langle \frac{y_{0}^{+}-y_{*}^+}{2},x_{N}^{+}\right\rangle \bigg],
			\end{aligned}
		\end{equation}
		where
		\begin{align*}
			\bfx & =\left[x_{0}^{+}-x_{0},x_{1}^{+}-x_{0}^{+},\ldots,x_{N}^{+}-x_{N-1}^{+}\right]^{\rmT},\\
			\bfy & =\left[y_{0}^{+}-y_{0},y_{1}^{+}-y_{0}^{+},\ldots,y_{N}^{+}-y_{N-1}^{+}\right]^{\rmT},
		\end{align*}
		is nonnegative.
	\end{enumerate}
\end{theorem}
\begin{proof}
	Note that the inequality \eqref{eq:F2-rate} holds for all [families $\{x_i\}_{i\in I_N^*}$ in $\mathbb{R}^d$, $\{g_i\}_{i\in I_N^*}$ in $\mathbb{R}^d$, and $\{f_i\}_{i\in I_N^*}$ in $\mathbb{R}$ satisfying $g_*=0$ and \eqref{eq:fsfom}] if and only if the inequality
	\begin{equation}
		\label{eq:appF2-2}
		\begin{aligned}
			& f_N^+ -f_*^+ + \frac{1}{2L}\Vert g_N\Vert^2-\frac{\tau}{2}\left\Vert x_{0}-x_{*}^+\right\Vert ^{2} \\
			& \qquad\leq\sum_{i,j\in I_N^*}\lambda_{ij}\left(f_{j}^{+}-f_{i}^{+}+\left\langle g_{j},x_{i}^{+}-x_{j}^{+}\right\rangle +\frac{\mu L}{2(L-\mu)}\left\Vert x_{i}^{+}-x_{j}^{+}\right\Vert ^{2}\right).
		\end{aligned}
	\end{equation}
	holds for all [families $\{x_0\}\cup\{x_i^+\}_{i\in I_N^*}$ in $\mathbb{R}^d$, $\{g_i\}_{i\in I_N^*}$ in $\mathbb{R}^d$, and $\{f_i^+\}_{i\in I_N^*}$ in $\mathbb{R}$ satisfying $g_*=0$ and \eqref{eq:fsfom-tilde}].\footnote{\label{footnote:F}To show ($\Rightarrow$), put $x_i=x_i^+ + \frac{1}{L}g_i$ and $f_i=f_i^+ + \frac{1}{2L}\|g_i\|^2$ for all $i\in I_N^*$. To show ($\Leftarrow$), put $x_i^+=x_i-\frac{1}{L}g_i$ and $f_i^+=f_i-\frac{1}{2L}\|g_i\|^2$ for all $i\in I_N^*$.} 
	
	For \eqref{eq:appF2-2} to hold for any $\{f_{k}^+\}_{k\in I_N^*}$, the terms $f_k^+$ in \eqref{eq:appF2-2} should vanish, which is \eqref{item:F2-1}. When \eqref{item:F2-1} holds, we can rewrite \eqref{eq:appF2-2} using $\bfg=-L\tilde{\bfH}^{-1}\bfx$ as
	\begin{align*}
		0&\leq S_1(x_0,x_0^+,\ldots,x_N^+,x_*^+)\\
		& :=\frac{\tau}{2}\left\Vert x_{0}-x_{*}^+\right\Vert ^{2}-\frac{L}{2}\left\Vert (\tilde{\bfH}^{-1}\bfx)_{N+1}\right\Vert ^{2}\\
		& \quad-L\left\langle \tilde{\bfH}^{-1}\bfx,\left[\begin{array}{c}
			\sum_{i\in I_N^*}\lambda_{i,0}\left(x_{i}^{+}-x_{0}^{+}\right)\\
			\sum_{i\in I_N^*}\lambda_{i,1}\left(x_{i}^{+}-x_{1}^{+}\right)\\
			\vdots\\
			\sum_{i\in I_N^*}\lambda_{i,N}\left(x_{i}^{+}-x_{N}^{+}\right)
		\end{array}\right]\right\rangle +\frac{\mu L}{2(L-\mu)}\sum_{i,j\in I_N^*}\lambda_{ij}\left\Vert x_{i}^{+}-x_{j}^{+}\right\Vert ^{2},
	\end{align*}
	where $(\tilde{\bfH}^{-1}\bfx)_{N+1}$ denotes the ($N+1$)-th entry of $\tilde{\bfH}^{-1}\bfx$. Thus, showing the inequality \eqref{eq:appF2-2} is equivalent to showing $\inf S_1\geq 0$, where the infimum is taken over all points $x_{0},x_{0}^{+},\ldots,x_{N}^{+},x_{*}^+$ in $\mathbb{R}^d$. We introduce \emph{dual} iterates $y_0^+,\ldots,y_N^+,y_*^+$ by the rule \eqref{eq:xmapy2}. Using \eqref{item:F2-1} and \eqref{eq:xmapy2}, we have
	\begin{equation}
		\label{eq:xmapy_lastrow}
		\begin{aligned}
			& \sum_{i\in I_N^*}\lambda_{i,N}\left(x_{i}^{+}-x_{N}^{+}\right) \\
			& \qquad=\sum_{i,j\in I_N^*}\lambda_{i,j}\left(x_{i}^{+}-x_{j}^{+}\right)-\sum_{j\in I_N^*\backslash\{N\}}\sum_{i\in I_N^*}\lambda_{i,j}\left(x_{i}^{+}-x_{j}^{+}\right)\\
			& \qquad=\left(x_{*}^+-x_{N}^{+}\right)-\sum_{k=1}^{N}\left(y_{k}^{+}-y_{k-1}^{+}\right)-\left(\frac{y_{0}^{+}+y_{*}^+}{2}-y_{N}^{+}-\frac{x_{N}^{+}-x_{*}^+}{2}\right)\\
			& \qquad=y_{0}^{+}-v,
		\end{aligned}
	\end{equation}
	where $v=\frac{y_{0}^{+}+y_{*}^+}{2}+\frac{x_{N}^{+}-x_{*}^+}{2}$. We can now rewrite \eqref{eq:appF2-2} as
	\begin{align*}
		0 & \leq S_{2}(x_{0},x_{0}^{+},\ldots,x_{N}^{+},x_{*}^+,y_{0}^{+},\ldots,y_{N}^{+},y_{*}^+)\\
		& :=\frac{\tau}{2}\left\Vert x_{0}-x_{*}^+\right\Vert ^{2}-\frac{L}{2}\left\Vert (\tilde{\bfH}^{-1}\bfx)_{N+1}\right\Vert ^{2}\\
		& \quad-L\left\langle \tilde{\bfH}^{-1}\bfx,\left[\begin{array}{c}
			y_{N}^{+}-y_{N-1}^{+}\\
			\vdots\\
			y_{1}^{+}-y_{0}^{+}\\
			y_{0}^{+}-v
		\end{array}\right]\right\rangle +\frac{\mu L}{2(L-\mu)}\sum_{i,j\in I_N^*}\lambda_{ij}\left\Vert x_{i}^{+}-x_{j}^{+}\right\Vert ^{2}.
	\end{align*}
	Thus, showing \eqref{eq:appF2-2} is equivalent to showing $\inf S_2\geq 0$, where the infimum is taken over all points $x_{0},x_{0}^{+},\ldots,x_{N}^{+},x_{*}^+,y_{0}^{+},\ldots,y_{N}^{+},y_{*}^+$ in $\mathbb{R}^d$ satisfying \eqref{eq:xmapy2}. Introduce an auxiliary point $y_0\in\mathbb{R}^d$ and consider the quantity
	\begin{equation}\nonumber
		\label{eq:appF2-4}
		\begin{aligned}
			& S_{3}(x_{0},x_{0}^{+},\ldots,x_{N}^{+},x_{*}^+,y_0,y_{0}^{+},\ldots,y_{N}^{+},y_{*}^+)  \\
			& \qquad:=\frac{L}{2}\left\Vert y_{0}-v\right\Vert ^{2}+\frac{\tau}{2}\left\Vert x_{0}-x_{*}^+\right\Vert ^{2}\\
			& \qquad\quad-L\left\langle \tilde{\bfH}^{-1}\bfx,\rmT\bfy\right\rangle +\frac{\mu L}{2(L-\mu)}\sum_{i,j\in I_N^*}\lambda_{ij}\left\Vert x_{j}^{+}-x_{i}^{+}\right\Vert ^{2}.
		\end{aligned}
	\end{equation}
	When the other points are fixed, one can verify that $S_2 = \inf_{y_0} S_3$. Thus, showing \eqref{eq:appF2-2} is equivalent to showing $\inf S_3 \geq 0$, where the infimum is taken over all points $x_{0},x_{0}^{+},\ldots,x_{N}^{+},x_{*}^+,y_0,y_{0}^{+},\ldots,y_{N}^{+},y_{*}^+$ in $\mathbb{R}^d$ satisfying \eqref{eq:xmapy2}. We proceed by rewriting the last term of $S_3$ as follows:
	\begin{align*}
		& \sum_{i,j\in I_N^*}\lambda_{i,j}\left\Vert x_{i}^{+}-x_{j}^{+}\right\Vert ^{2}\\
		& \qquad=\sum_{k\in I_N^*}\left(\sum_{j\in I_N^*}\lambda_{k,j}+\sum_{i\in I_N^*}\lambda_{i,k}\right)\left\Vert x_{k}^{+}\right\Vert ^{2}-2\sum_{j\in I_N^*}\sum_{i\in I_N^*}\lambda_{i,j}\left\langle x_{i}^{+},x_{j}^{+}\right\rangle \\
		& \qquad=\sum_{k\in I_N^*}\left(\sum_{j\in I_N^*}\lambda_{k,j}-\sum_{i\in I_N^*}\lambda_{i,k}\right)\left\Vert x_{k}^{+}\right\Vert ^{2}-2\sum_{j\in I_N^*}\left\langle \sum_{i\in I_N^*}\lambda_{i,j}\left(x_{i}^{+}-x_{j}^{+}\right),x_{j}^{+}\right\rangle .
	\end{align*}
	We simplify the first term using \eqref{item:F2-1}, and the second term using \eqref{eq:xmapy2} and \eqref{eq:xmapy_lastrow}.
	\begin{align*}
		& \sum_{i,j\in I_N^*}\lambda_{i,j}\left\Vert x_{i}^{+}-x_{j}^{+}\right\Vert ^{2} \\
		& \qquad=-\left\Vert x_{N}^{+}\right\Vert ^{2}+\left\Vert x_{*}^+\right\Vert ^{2}-2\sum_{k=1}^{N}\left\langle y_{k}^{+}-y_{k-1}^{+},x_{N-k}^{+}\right\rangle \\
		& \qquad\quad-2\left\langle \frac{y_{0}^{+}+y_{*}^+}{2}-y_{N}^{+}-\frac{x_{N}^{+}-x_{*}^+}{2},x_{*}^+\right\rangle -2\left\langle \frac{y_{0}^{+}-y_{*}^+}{2}-\frac{x_{N}^{+}-x_{*}^+}{2},x_{N}^{+}\right\rangle \\
		& \qquad=-2\left(\sum_{k=1}^{N}\left\langle y_{k}^{+}-y_{k-1}^{+},x_{N-k}^{+}\right\rangle +\left\langle \frac{y_{0}^{+}+y_{*}^+}{2}-y_{N}^{+},x_{*}^+\right\rangle +\left\langle \frac{y_{0}^{+}-y_{*}^+}{2},x_{N}^{+}\right\rangle \right).
	\end{align*}
	This completes the proof.
\end{proof}

\subsection{PEP analysis for performance criterion ${\|\nabla f(x_{N})\|^{2}} / {(f(x_0)-f_*)}$}
\label{sec:G2}

In this subsection, we present a PEP analysis aimed at proving convergence guarantees in the form $\frac{1}{2\tau}\|\nabla f(y_{N})\|^{2}\leq f(y_{0})-f_*$. Our analysis differs from the one in \cite{kim2021optimizing} by using $\{y_{k+1}^+ - y_k^+\}$ instead of $\{g_k\}$ as the basis for the quadratic form (see Footnote~\ref{footnote:see-psd}), although both analyses aim at the same performance criterion.
\begin{theorem}
	\label{thm:G2}
	Given a set of nonnegative multipliers $\{\nu_{i,j}\}_{i,j\in I_N^*}$, the inequality
	\begin{equation}
		\label{eq:G2-rate}
		\begin{aligned}
			& \frac{1}{2\tau}\left\Vert g_{N}\right\Vert ^{2}-f_{0}+f_{*}-\frac{L}{2}\left\Vert \sum_{i\in I_{N}^{*}}\nu_{i,*}\left(y_{i}-\frac{1}{L}g_{i}-y_{*}\right)\right\Vert ^{2}\\
			& \qquad\leq\sum_{i,j\in I_{N}^{*}}\nu_{i,j}\bigg[f_{j}-f_{i}+\left\langle g_{j},y_{i}-y_{j}\right\rangle +\frac{1}{2L}\left\Vert g_{i}-g_{j}\right\Vert ^{2}\\
			& \qquad\qquad\qquad\qquad\quad+\frac{\mu L}{2(L-\mu)}\left\Vert y_{i}-y_{j}-\frac{1}{L}\left(g_{i}-g_{j}\right)\right\Vert ^{2}\bigg]
		\end{aligned}
	\end{equation}
	holds for all [families $\{y_i\}_{i\in I_N^*}$ in $\mathbb{R}^d$, $\{g_i\}_{i\in I_N^*}$ in $\mathbb{R}^d$, and $\{f_i\}_{i\in I_N^*}$ in $\mathbb{R}$ satisfying $g_*=0$ and \eqref{eq:fsfom-dual}] if and only if the following two conditions are satisfied:
	\begin{enumerate}[label=\textnormal{(G\arabic*)},ref=G\arabic*]
		\item \label{item:G2-1} \[
		\sum_{i\in I_N^*}\nu_{i,k}-\sum_{j\in I_N^*}\nu_{k,j}=\begin{cases}
			1 & \textup{if }k=*,\\
			-1 & \textup{if }k=0,\\
			0 & \textup{otherwise}.
		\end{cases}
		\]
		\item \label{item:G2-2} For any points $x_{0},x_{0}^{+},\ldots,x_{N}^{+},x_*^+,y_{0},y_{0}^{+},\ldots,y_{N}^{+},y_*^+\in\mathbb{R}^d$ satisfying
		\begin{equation}
			\label{eq:ymapx2}
			\begin{aligned}
				x_{k}^{+}-x_{k-1}^{+}& =\sum_{i\in I_N^*}\nu_{i,N-k}\left(y_{i}^{+}-y_{N-k}^{+}\right) \quad\textup{for }k=1,\ldots,N-1,\\
				x_{0}^{+}-x_{*}^+& =\sum_{i\in I_N^*}\nu_{i,N}\left(y_{i}^{+}-y_{N}^{+}\right),\\
				\frac{x_{*}^+-x_{N}^{+}}{2}& =\frac{y_{*}-y_{0}^{+}}{2}+\sum_{i\in I_N^*}\nu_{i,*}\left(y_{i}^{+}-y_{*}^+\right),
			\end{aligned}
		\end{equation}
		the quantity \eqref{eq:F-PEP-functional2} is nonnegative.
	\end{enumerate}
\end{theorem}
\begin{proof}
	Note that the inequality \eqref{eq:G2-rate} holds for all [families $\{y_i\}_{i\in I_N^*}$ in $\mathbb{R}^d$, $\{g_i\}_{i\in I_N^*}$ in $\mathbb{R}^d$, and $\{f_i\}_{i\in I_N^*}$ in $\mathbb{R}$ satisfying $g_*=0$ and \eqref{eq:fsfom-dual}] if and only if the inequality
	\begin{equation}
		\label{eq:appG2-2}
		\begin{aligned}
			& \frac{1}{2\tau}\left\Vert g_{N}\right\Vert ^{2}+f_{*}^{+}-f_{0}^{+}-\frac{1}{2L}\Vert g_{0}\Vert^{2}-\frac{L}{2}\left\Vert \sum_{i\in I_{N}^{*}}\nu_{i,*}\left(y_{i}^{+}-y_{*}^{+}\right)\right\Vert ^{2}\\
			& \qquad\leq\sum_{i,j\in I_{N}^{*}}\nu_{ij}\left(f_{j}^{+}-f_{i}^{+}+\left\langle g_{j},y_{i}^{+}-y_{j}^{+}\right\rangle +\frac{\mu L}{2(L-\mu)}\left\Vert y_{i}^{+}-y_{j}^{+}\right\Vert ^{2}\right)
		\end{aligned}
	\end{equation}
	holds for all [families $\{y_0\}\cup\{y_i^+\}_{i\in I_N^*}$ in $\mathbb{R}^d$, $\{g_i\}_{i\in I_N^*}$ in $\mathbb{R}^d$, and $\{f_i^+\}_{i\in I_N^*}$ in $\mathbb{R}$ satisfying $g_*=0$ and \eqref{eq:fsfom-dual-tilde}] (see Footnote~\ref{footnote:F}).
	
	For \eqref{eq:appG2-2} to hold for any $\{f_{k}^+\}_{k\in I_N^*}$, the terms $f_k^+$ in \eqref{eq:appG2-2} should vanish, which is \eqref{item:G2-1}. When \eqref{item:G2-1} holds, we can rewrite \eqref{eq:appG2-2} using $\bfg=-L\tilde{\bfH}^{-\rmA}\bfy$ (in particular, we have $g_{0}=-L(y_{0}^{+}-y_{0})$ since $\tilde{\bfH}^{-\rmA}_{1,1}=1$) as
	\begin{align*}
		0 & \leq S_{1}(y_{0},y_{0}^{+},\ldots,y_{N}^{+},y_{*}^+)\\
		& :=\frac{L}{2}\left\Vert y_0^+-y_0\right\Vert ^{2}-\frac{L^{2}}{2\tau}\left\Vert (\tilde{\bfH}^{-\rmA}\bfy)_{N+1}\right\Vert ^{2}+\frac{L}{2}\left\Vert \sum_{i\in I_N^*}\nu_{i,*}\left(y_{i}^{+}-y_{*}^+\right)\right\Vert ^{2}\\
		& \quad-L\left\langle \tilde{\bfH}^{-\rmA}\bfy,\left[\begin{array}{c}
			\sum_{i\in I_N^*}\nu_{i,0}\left(y_{i}^{+}-y_{0}^{+}\right)\\
			\sum_{i\in I_N^*}\nu_{i,1}\left(y_{i}^{+}-y_{1}^{+}\right)\\
			\vdots\\
			\sum_{i\in I_N^*}\nu_{i,N}\left(y_{i}^{+}-y_{N}^{+}\right)
		\end{array}\right]\right\rangle +\frac{\mu L}{2(L-\mu)}\sum_{i,j\in I_N^*}\nu_{ij}\left\Vert y_{i}^{+}-y_{j}^{+}\right\Vert ^{2},
	\end{align*}
	where $(\tilde{\bfH}^{-\rmA}\bfy)_{N+1}$ denotes the ($N+1$)-th entry of $\tilde{\bfH}^{-\rmA}\bfy$. Thus, showing the inequality \eqref{eq:appG2-2} is equivalent to showing $\inf S_1\geq 0$, where the infimum is taken over all points $y_0,y_0^+,\ldots,y_N^+,y_*^+$ in $\mathbb{R}^d$. We introduce dual iterates $x_0^+,\ldots,x_N^+,x_*^+$ by the rule \eqref{eq:ymapx2}. Using \eqref{item:G2-1} and \eqref{eq:ymapx2}, we have
	\begin{equation}
		\label{eq:ymapx_lastrow}
		\begin{aligned}
			\sum_{i\in I_{N}^{*}}\nu_{i,0}\left(y_{i}^{+}-y_{0}^{+}\right) & =\sum_{i,j\in I_{N}^{*}}\nu_{i,j}\left(y_{i}^{+}-y_{j}^{+}\right)-\sum_{j\in I_{N}^{*}\backslash\{0\}}\sum_{i\in I_{N}^{*}}\nu_{i,j}\left(y_{i}^{+}-y_{j}^{+}\right)\\
			& =\left(y_{0}^{+}-y_{*}^{+}\right)-\sum_{k=1}^{N-1}\left(x_{k}^{+}-x_{k-1}^{+}\right)\\
			& \quad-\left(x_{0}^{+}-x_{*}^{+}\right)-\left(\frac{x_{*}^{+}-x_{N}^{+}}{2}-\frac{y_{*}^{+}-y_{0}^{+}}{2}\right)\\
			& =x_{N}^{+}-x_{N-1}^{+}+y_{0}^{+}-v,
		\end{aligned}
	\end{equation}
	where $v=\frac{y_{0}^{+}+y_{*}^+}{2}+\frac{x_{N}^{+}-x_{*}^+}{2}$. Note that the last equality in \eqref{eq:ymapx2} can be written as $\sum_{i\in I_N^*}\nu_{i,*}(y_{i}^{+}-y_{*}^{+})=y_{0}^+-v$. We can now rewrite \eqref{eq:appG2-2} as
	\begin{align*}
		0 & \leq S_{2}(x_{0}^{+},\ldots,x_{N}^{+},x_{*},y_{0},y_{0}^{+},\ldots,y_{N}^{+},y_{*})\\
		& :=\frac{L}{2}\left\Vert y_{0}^{+}-y_{0}\right\Vert ^{2}-\frac{L^{2}}{2\tau}\left\Vert (\tilde{\bfH}^{-\rmA}\bfy)_{N+1}\right\Vert ^{2}+\frac{L}{2}\left\Vert y_{0}^{+}-v\right\Vert ^{2}\\
		& \quad-L\left\langle \tilde{\bfH}^{-\rmA}\bfy,\left[\begin{array}{c}
			x_{N}^{+}-x_{N-1}^{+}+y_{0}^{+}-v\\
			\vdots\\
			x_{1}^{+}-x_{0}^{+}\\
			x_{0}^{+}-x_{*}^+
		\end{array}\right]\right\rangle +\frac{\mu L}{2(L-\mu)}\sum_{i,j\in I_N^*}\nu_{ij}\left\Vert y_{i}^{+}-y_{j}^{+}\right\Vert ^{2}\\
		& =\frac{L}{2}\left\Vert y_{0}-v\right\Vert ^{2}-\frac{L^{2}}{2\tau}\left\Vert (\tilde{\bfH}^{-\rmA}\bfy)_{N+1}\right\Vert ^{2}\\
		& \quad-L\left\langle \tilde{\bfH}^{-\rmA}\bfy,\left[\begin{array}{c}
			x_{N}^{+}-x_{N-1}^{+}\\
			\vdots\\
			x_{1}^{+}-x_{0}^{+}\\
			x_{0}^{+}-x_{*}^+
		\end{array}\right]\right\rangle +\frac{\mu L}{2(L-\mu)}\sum_{i,j\in I_N^*}\nu_{ij}\left\Vert y_{i}^{+}-y_{j}^{+}\right\Vert ^{2}.
	\end{align*}
	Thus, showing \eqref{eq:appG2-2} is equivalent to showing $\inf S_2\geq 0$, where the infimum is taken over all points $x_{0}^{+},\ldots,x_{N}^{+},x_{*}^+,y_0,y_{0}^{+},\ldots,y_{N}^{+},y_{*}^+$ satisfying \eqref{eq:ymapx2}. Introduce an auxiliary point $x_0\in\mathbb{R}^d$ and consider the quantity
	\begin{align*}
		& S_{3}(x_{0},x_{0}^{+},\ldots,x_{N}^{+},x_{*}^+,y_{0},y_{0}^{+},\ldots,y_{N}^{+},y_{*}^+)\\
		& \qquad:=\frac{L}{2}\left\Vert y_{0}-v\right\Vert ^{2}+\frac{\tau}{2}\left\Vert x_{0}-x_{*}^+\right\Vert ^{2}\\
		& \qquad\quad-L\left\langle \tilde{\bfH}^{-\rmA}\bfy,\rmT\bfx\right\rangle +\frac{\mu L}{2(L-\mu)}\sum_{i,j\in I_N^*}\nu_{ij}\left\Vert y_{i}^{+}-y_{j}^{+}\right\Vert ^{2}.
	\end{align*}
	When the other points are fixed, one can verify that $S_2 = \inf_{x_0} S_3$. Thus, showing \eqref{eq:appG2-2} is equivalent to showing $S_3 \geq 0$, where the infimum is taken over all points $x_{0},x_{0}^{+},\ldots,x_{N}^{+},x_{*},y_0,y_{0}^{+},\ldots,y_{N}^{+},y_{*}$ satisfying \eqref{eq:ymapx2}. By Definition~\ref{def:at}, we have $\langle \tilde{\bfH}^{-\rmA}\bfy,\rmT\bfx\rangle =\langle \rmT\bfy,\tilde{\bfH}^{-1}\bfx\rangle $. We proceed by rewriting the last term of $S_3$ as follows:
	\begin{align*}
		& \sum_{i,j\in I_N^*}\nu_{i,j}\left\Vert y_{i}^{+}-y_{j}^{+}\right\Vert ^{2}\\
		&\qquad =\sum_{k\in I_N^*}\left(\sum_{j\in I_N^*}\nu_{k,j}+\sum_{i\in I_N^*}\nu_{i,k}\right)\left\Vert y_{k}^{+}\right\Vert ^{2}-2\sum_{j\in I_N^*}\sum_{i\in I_N^*}\nu_{i,j}\left\langle y_{i}^{+},y_{j}^{+}\right\rangle \\
		&\qquad =\sum_{k\in I_N^*}\left(\sum_{j\in I_N^*}\nu_{k,j}-\sum_{i\in I_N^*}\nu_{i,k}\right)\left\Vert y_{k}^{+}\right\Vert ^{2}-2\sum_{j\in I_N^*}\left\langle \sum_{i\in I_N^*}\nu_{i,j}\left(y_{i}^{+}-y_{j}^{+}\right),y_{j}^{+}\right\rangle .
	\end{align*}
	We simplify the first term using \eqref{item:G2-1}, and the second term using \eqref{eq:ymapx2} and \eqref{eq:ymapx_lastrow}.
	\begin{align*}
		& \sum_{i,j\in I_N^*}\nu_{i,j}\left\Vert y_{i}^{+}-y_{j}^{+}\right\Vert ^{2} \\
		&\qquad =-\left\Vert y_{*}^+\right\Vert ^{2}+\left\Vert y_{0}^{+}\right\Vert ^{2}-2\sum_{k=1}^{N-1}\left\langle x_{k}^{+}-x_{k-1}^{+},y_{N-k}^{+}\right\rangle -2\left\langle x_{0}^{+}-x_{*}^+,y_{N}^{+}\right\rangle \\
		&\qquad \quad-2\left\langle \frac{x_{*}^+-x_{N}^{+}}{2}-\frac{y_{*}-y_{0}^{+}}{2},y_{*}^+\right\rangle -2\left\langle \frac{x_{N}^{+}+x_{*}^+}{2}-x_{N-1}^{+}+\frac{y_{0}^{+}-y_{*}^+}{2},y_{0}^{+}\right\rangle \\
		&\qquad =-2\left(\sum_{k=1}^{N}\left\langle y_{k}^{+}-y_{k-1}^{+},x_{N-k}^{+}\right\rangle +\left\langle \frac{y_{0}^{+}+y_{*}^+}{2}-y_{N}^{+},x_{*}^+\right\rangle +\left\langle \frac{y_{0}^{+}-y_{*}^+}{2},x_{N}^{+}\right\rangle \right).
	\end{align*}
	This completes the proof.
\end{proof}
The presence of the last term on the left-hand side of \eqref{eq:G2-rate} is undesirable. The following theorem is a variant of Theorem~\ref{thm:G2}, where this term is removed under an additional assumption that the point $x_*$ is involved in the convergence analysis only through the interpolation inequality \eqref{eq:interpolation_ineqq} with $p=x_*$ and $q=x_N$. 
\begin{theorem}
	\label{thm:G2-modified}
	Given a set of nonnegative multipliers $\{\nu_{i,j}\}_{i,j\in I_N^*}$ satisfying
	\begin{equation}
		\label{eq:simple-g-condition}
		\nu_{i,*}=\begin{cases}
			1 & \textup{if }i=N,\\
			0 & \textup{otherwise},
		\end{cases}\qquad\nu_{*,j}=0\quad\textup{for } j=0,\ldots N,
	\end{equation}
	the inequality
	\begin{equation}
		\label{eq:G2-rate-modified}
		\begin{aligned}
			\frac{1}{2\tau}\left\Vert g_{N}\right\Vert ^{2}-f_{0}+f_{*} & \leq\sum_{i,j\in I_{N}^{*}}\nu_{i,j}\bigg[f_{j}-f_{i}+\left\langle g_{j},y_{i}-y_{j}\right\rangle +\frac{1}{2L}\left\Vert g_{i}-g_{j}\right\Vert ^{2}\\
			& \qquad\qquad\qquad\quad+\frac{\mu L}{2(L-\mu)}\left\Vert y_{i}-y_{j}-\frac{1}{L}\left(g_{i}-g_{j}\right)\right\Vert ^{2}\bigg]
		\end{aligned}
	\end{equation}
	holds for all [families $\{y_i\}_{i\in I_N^*}$ in $\mathbb{R}^d$, $\{g_i\}_{i\in I_N^*}$ in $\mathbb{R}^d$, and $\{f_i\}_{i\in I_N^*}$ in $\mathbb{R}$ satisfying $g_*=0$ and \eqref{eq:fsfom-dual}] if and only if the two conditions \eqref{item:G2-1} and \eqref{item:G2-2} are satisfied.
\end{theorem}
\begin{proof}
	Assume that \eqref{eq:simple-g-condition} holds. Then, the quantity $S_1$ in the original proof of Theorem~\ref{thm:G2} is the sum of $(\frac{L}{2}+\frac{\mu L}{2(L-\mu)})\Vert y_{*}^+-y_{N}^{+}\Vert ^{2}$ and other terms independent of $y_*^+$, and is thus minimized at $y_*^+=y_N^+$ when the other points $y_0,y_0^+,\ldots,y_N^+$ are fixed. 
	
	If the inequality \eqref{eq:G2-rate-modified} is used instead of \eqref{eq:G2-rate}, mimicking the proof of Theorem~\ref{thm:G2} yields the following quantity instead of $S_1$:
	\begin{align*}
		S_{1}'(y_{0},y_{0}^{+},\ldots,y_{N}^{+},y_{*}^{+}) & :=\frac{L}{2}\left\Vert y_{0}^{+}-y_{0}\right\Vert ^{2}-\frac{L^{2}}{2\tau}\left\Vert (\tilde{\bfH}^{-\rmA}\bfy)_{N+1}\right\Vert ^{2}\\
		& \quad-L\left\langle \tilde{\bfH}^{-\rmA}\bfy,\left[\begin{array}{c}
			\sum_{i\in I_{N}^{*}}\nu_{i,0}\left(y_{i}^{+}-y_{0}^{+}\right)\\
			\sum_{i\in I_{N}^{*}}\nu_{i,1}\left(y_{i}^{+}-y_{1}^{+}\right)\\
			\vdots\\
			\sum_{i\in I_{N}^{*}}\nu_{i,N}\left(y_{i}^{+}-y_{N}^{+}\right)
		\end{array}\right]\right\rangle \\
		& \quad+\frac{\mu L}{2(L-\mu)}\sum_{i,j\in I_{N}^{*}}\nu_{ij}\left\Vert y_{i}^{+}-y_{j}^{+}\right\Vert ^{2}.
	\end{align*}
	Since $S_1'$ is the sum of $\frac{\mu L}{2(L-\mu)}\Vert y_{*}^+-y_{N}^{+}\Vert ^{2}$ and other terms independent of $y_*^+$, it is minimized at $y_*^+=y_N^+$ when the other points $y_0,y_0^+,\ldots,y_N^+$ are fixed. 
	
	Thus, we have
	\begin{align*}
		\inf_{y_{0},y_{0}^{+},\ldots,y_{N}^{+},y_{*}^+\in\mathbb{R}^d}S_{1}(y_{0},y_{0}^{+},\ldots,y_{N}^{+},y_{*}^+) & =\inf_{y_{0},y_{0}^{+},\ldots,y_{N}^{+}\in\mathbb{R}^d}S_{1}(y_{0},y_{0}^{+},\ldots,y_{N}^{+},y_{N}^{+})\\
		& =\inf_{y_{0},y_{0}^{+},\ldots,y_{N}^{+}\in\mathbb{R}^d}S_{1}'(y_{0},y_{0}^{+},\ldots,y_{N}^{+},y_{N}^{+})\\
		& =\inf_{y_{0},y_{0}^{+},\ldots,y_{N}^{+},y_{*}^+\in\mathbb{R}^d}S_{1}'(y_{0},y_{0}^{+},\ldots,y_{N}^{+},y_{*}^+).
	\end{align*}
	Therefore, showing $\inf S_1'\geq 0$ is equivalent to showing $\inf S_1\geq 0$. This means that Theorem~\ref{thm:G2} remains valid after replacing the inequality \eqref{eq:G2-rate} by \eqref{eq:G2-rate-modified}.
\end{proof}

\subsection{General correspondence between two PEP analyses}
\label{sec:cor2}

In this subsection, we present a technique to convert a PEP analysis using Theorem~\ref{thm:G2} (or Theorem~\ref{thm:G2-modified}) into a PEP analysis using Theorem~\ref{thm:F2}, and vice versa. Suppose that we have a set of multipliers $\{\nu_{i,j}\}_{i,j\in I_N^*}$ satisfying the two conditions \eqref{item:G2-1} and \eqref{item:G2-2}, and that we want to find a set of nonnegative multipliers $\{\lambda_{i,j}\}_{i,j\in I_N^*}$ satisfying the two conditions \eqref{item:F2-1} and \eqref{item:F2-2}. Then, we take the following two steps:

{\it Step} 1. For $j\neq N$, set the multipliers $\lambda_{i,j}$ so that \eqref{eq:xmapy2} and \eqref{eq:ymapx2} are equivalent.

{\it Step} 2. For $j=N$, set the multipliers $\lambda_{i,j}$ so that \eqref{item:F2-1} holds. 

If these two steps are completed successfully, we can conclude that the two conditions \eqref{item:F2-1} and \eqref{item:F2-2} are satisfied for the resulting set of multipliers $\{\lambda_{i,j}\}_{i,j\in I_N^*}$. 

For the opposite direction, we take the following two steps:

{\it Step} 1. For $j\neq 0$, set the multipliers $\nu_{i,j}$ so that \eqref{eq:xmapy2} and \eqref{eq:ymapx2} are equivalent.

{\it Step} 2. For $j=0$, set the multipliers $\nu_{i,j}$ so that \eqref{item:G2-1} holds.

\begin{remark}
	For a given set of non-negative multipliers $\{\nu_{i,j}\}_{i,j\in I_N^*}$, there might not exist a set of non-negative multipliers $\{\lambda_{i,j}\}_{i,j\in I_N^*}$ for which \eqref{eq:xmapy2} and \eqref{eq:ymapx2} are equivalent, and vice versa. Thus, not all PEP analyses using Theorem~\ref{thm:G2} (or Theorem~\ref{thm:G2-modified}) can be translated to PEP analyses using Theorem~\ref{thm:F2}, and vice versa.
\end{remark}

\subsection{Special correspondence between two relaxed PEP analyses}
\label{sec:restricted}

The results in this subsection will not be directly used in the rest of the paper, thus the reader can skip it. In this subsection, we show a full correspondence between every valid PEP analysis using Theorem~\ref{thm:F2} and every valid PEP analysis using Theorem~\ref{thm:G2-modified}, if it is further assumed that only a specific subset of interpolation inequalities can be used. Specifically, we \emph{relax} the analysis in Section~\ref{sec:F2} by assuming \eqref{eq:relaxed-f}, and the analysis in Section~\ref{sec:G2} by assuming
\begin{equation}
	\label{eq:relaxed-g}
	\nu_{i,j}=0\quad\textup{for all }(i,j)\notin\left\{ (k,k+1),(N,k)\right\} _{0\leq k\leq N-1}\cup\left\{ (N,*)\right\} .
\end{equation}
By applying the technique in the previous subsection, we obtain the following lemma.

\begin{lemma}
	\label{lem:relaxed}\ 
	\begin{enumerate}[label=\textnormal{(\roman*)}]
		\item \label{lem:relaxed-item1} If a set of nonnegative multipliers $\{\nu_{i,j}\}_{i,j\in I_N^*}$ satisfies \eqref{item:G2-1} and \eqref{eq:relaxed-g}, then there is a set of nonnegative multipliers $\{\lambda_{i,j}\}_{i,j\in I_N^*}$ satisfying \eqref{item:F2-1} and \eqref{eq:relaxed-f}, such that \eqref{eq:xmapy2} and \eqref{eq:ymapx2} are equivalent.
		\item \label{lem:relaxed-item2} If a set of nonnegative multipliers $\{\lambda_{i,j}\}_{i,j\in I_N^*}$ satisfies \eqref{item:F2-1} and \eqref{eq:relaxed-f}, then there is a set of nonnegative multipliers $\{\nu_{i,j}\}_{i,j\in I_N^*}$ satisfying \eqref{item:G2-1} and \eqref{eq:relaxed-g}, such that \eqref{eq:xmapy2} and \eqref{eq:ymapx2} are equivalent.
	\end{enumerate}
\end{lemma}
\begin{proof}
	(i) Let $A_{k}=1+\sum_{j=0}^{k-1}\nu_{N,j}$ for $k=1,\ldots,N$. Then, $\{A_k\}_{k=1}^N$ is a nondecreasing sequence in $[1,\infty)$. Using \eqref{eq:relaxed-g} and \eqref{item:G2-1}, we can express all the multipliers $\nu_{i,j}$ in terms of $A_k$ as
	\begin{equation}
		\label{eq:A-to-nu}
		\nu_{i,j}=\begin{cases}
			A_{j} & \textup{if }i=j-1,\ 1\leq j\leq N,\\
			A_{1}-1 & \textup{if }i=N,\ j=0,\\
			A_{j+1}-A_{j} & \textup{if }i=N,\ 1\leq j\leq N-1,\\
			1 & \textup{if }i=N,\ j=*,\\
			0 & \textup{otherwise}.
		\end{cases}
	\end{equation}
	Thus, \eqref{eq:ymapx2} without the last equality can be written as
	\begin{equation}
		\label{eq:transform_matrix}
		\left[\begin{array}{c}
			x_{N-1}^{+}-x_{N-2}^{+}\\
			\vdots\\
			x_{1}^{+}-x_{0}^{+}\\
			x_{0}^{+}-x_{*}^+
		\end{array}\right]=-\left[\begin{array}{ccc}
			\bfM_{1,1} & \cdots & \bfM_{1,N}\\
			& \ddots & \vdots\\
			&  & \bfM_{N,N}
		\end{array}\right]\left[\begin{array}{c}
			y_{1}^{+}-y_{0}^{+}\\
			\vdots\\
			y_{N-1}^{+}-y_{N-2}^{+}\\
			y_{N}^{+}-y_{N-1}^{+}
		\end{array}\right],
	\end{equation}
	where
	\begin{equation}
		\label{eq:M-choice}
		\bfM_{i,j}=\begin{cases}
			A_{i} & \textup{if }j=i,\\
			A_{i}-A_{i+1} & \textup{if }j\geq i+1.
		\end{cases}
	\end{equation}
	The last equality of \eqref{eq:ymapx2} can be written as
	\begin{equation}
		\label{eq:last_simple}
		\frac{x_{*}^+-x_{N}^{+}}{2}=y_{N}^{+}-\frac{y_{*}^++y_{0}^{+}}{2}.
	\end{equation}
	We can rewrite \eqref{eq:transform_matrix} using the inverse matrix $\bfM^{-1}$ as
	\begin{equation}
		\label{eq:transform-inverse}
		\left[\begin{array}{c}
			y_{1}^{+}-y_{0}^{+}\\
			\vdots\\
			y_{N-1}^{+}-y_{N-2}^{+}\\
			y_{N}^{+}-y_{N-1}^{+}
		\end{array}\right]=-\left[\begin{array}{ccc}
			\bfM_{1,1}^{-1} & \cdots & \bfM_{1,N}^{-1}\\
			& \ddots & \vdots\\
			&  & \bfM_{N,N}^{-1}
		\end{array}\right]\left[\begin{array}{c}
			x_{N-1}^{+}-x_{N-2}^{+}\\
			\vdots\\
			x_{1}^{+}-x_{0}^{+}\\
			x_{0}^{+}-x_{*}
		\end{array}\right],
	\end{equation}
	where
	\begin{equation}
		\label{eq:Minv-choice}
		\bfM_{i,j}^{-1}=\begin{cases}
			{1}/{A_{i}} & \textup{if }j=i,\\
			{1}/{A_{i}}-1/{A_{i+1}} & \textup{if }j\geq i+1.
		\end{cases}
	\end{equation}
	We set the multipliers $\lambda_{i,j}$ as
	\begin{equation}
		\label{eq:choice-of-lam-simple}
		\lambda_{i,j}=\begin{cases}
			1/A_{N-j+1} & \textup{if }i=j-1,\ 1\leq j\leq N,\\
			1/A_{N} & \textup{if }i=*,\ j=0,\\
			1/A_{N-j}-1/A_{N-j+1} & \textup{if }i=*,\ 1\leq j\leq N-1,\\
			1-1/A_{1} & \textup{if }i=*,\ j=N,\\
			0 & \textup{otherwise}.
		\end{cases}
	\end{equation}
	Then, the condition \eqref{eq:xmapy2} is equivalent to \eqref{eq:last_simple} and \eqref{eq:transform-inverse} with \eqref{eq:Minv-choice}. Thus, \eqref{eq:xmapy2} is equivalent to \eqref{eq:ymapx2}. It is straightforward to check \eqref{item:F2-1}, \eqref{eq:relaxed-f}, and that $\lambda_{i,j}\geq 0$ for all $i,j\in I_N^*$.
	
	(ii) Let $A_{k}=1/\sum_{j=0}^{N-k} \lambda_{*,j}$ for $k=1,\ldots, N$. Then, $\{A_k\}_{k=1}^{N}$ is a nondecreasing sequence in $[1,\infty)$. Using \eqref{eq:relaxed-f} and \eqref{item:F2-1}, we can express all the multipliers $\lambda_{i,j}$ in terms of $A_k$ as in \eqref{eq:choice-of-lam-simple}. We set the multipliers $\nu_{i,j}$ as in \eqref{eq:A-to-nu}. Reversing the argument of the proof for the case (i), one can see that this set of multipliers $\{\lambda_{i,j}\}_{i,j\in I_N^*}$ satisfies the desired properties.
\end{proof}

The following theorem follows from Lemma~\ref{lem:relaxed}, Theorems~\ref{thm:F2}, and \ref{thm:G2-modified}.
\begin{theorem}
	\label{thm:relaxed}
	The following statements are equivalent:
	\begin{enumerate}[label=\textnormal{(\roman*)}]
		\item There is a set of nonnegative multipliers $\{\nu_{i,j}\}_{i,j\in I_N^*}$ satisfying \eqref{eq:relaxed-g}, such that the inequality \eqref{eq:G2-rate-modified} holds for all [families $\{y_i\}_{i\in I_N^*}$ in $\mathbb{R}^d$, $\{g_i\}_{i\in I_N^*}$ in $\mathbb{R}^d$, and $\{f_i\}_{i\in I_N^*}$ in $\mathbb{R}$ satisfying $g_*=0$ and \eqref{eq:fsfom-dual}].
		\item There is a set of nonnegative multipliers $\{\lambda_{i,j}\}_{i,j\in I_N^*}$ satisfying \eqref{eq:relaxed-f}, such that the inequality \eqref{eq:F2-rate} holds for all [families $\{x_i\}_{i\in I_N^*}$ in $\mathbb{R}^d$, $\{g_i\}_{i\in I_N^*}$ in $\mathbb{R}^d$, and $\{f_i\}_{i\in I_N^*}$ in $\mathbb{R}$ satisfying $g_*=0$ and \eqref{eq:fsfom}].
	\end{enumerate}
\end{theorem}
Theorem~\ref{thm:relaxed} is valid for any $\mu \in [0, L)$. See \cite[Thm.~1]{kim2024time} for the same conclusion in the nonstrongly convex case ($\mu=0$). See \cite[Thm.~3]{kim2024convergence} and \cite[Thm.~2]{kim2024time} for the essentially same conclusion in the corresponding continuous-time analysis ($\mu=0$ and $1/L\to0$). 

\subsection{Special correspondence that will connect Theorems~\ref{conj:gd} and \ref{thm:gd}}
\label{sec:restricted-2}

In this subsection, we present a variant of the (i) $\Rightarrow$ (ii) direction of Theorem~\ref{thm:relaxed} with a weaker assumption and a weaker conclusion. Precisely, we replace the condition \eqref{eq:relaxed-g} with the following one:
\begin{equation}
	\label{eq:d1-tau}
	\nu_{i,j}= 0\quad\textup{for all }(i,j)\notin\left\{ (k,k+1),(k+1,k),(N,k)\right\} _{0\leq k\leq N-1}\cup\left\{ (N,*)\right\} .
\end{equation}
The following lemma is a variant of Lemma~\ref{lem:relaxed}\ref{lem:relaxed-item1} in this context.
\begin{lemma}
	\label{lem:relaxed_general}
	If a set of nonnegative multipliers $\{\nu_{i,j}\}_{i,j\in I_N^*}$ satisfies \eqref{item:G2-1} and \eqref{eq:d1-tau}, then there is a set of nonnegative multipliers $\{\lambda_{i,j}\}_{i,j\in I_N^*}$ satisfying \eqref{item:F2-1} such that \eqref{eq:xmapy2} is equivalent to \eqref{eq:ymapx2}.
\end{lemma}
The proof of Lemma~\ref{lem:relaxed_general} follows a similar argument as in the proof of Lemma~\ref{lem:relaxed}, but involves a tedious computation. In particular, we explicitly compute the multipliers $\lambda_{i,j}$ in terms of $\nu_{i,j}$. Thus, we defer the proof to Appendix~\ref{app:lem:relaxed_general}. By combining Lemma~\ref{lem:relaxed_general} with Theorems~\ref{thm:F2} and \ref{thm:G2-modified}, we obtain the following result.
\begin{theorem}
	\label{thm:relaxed_general}
	We have \textup{(i)} $\Rightarrow$ \textup{(ii)} for the following statements:
	\begin{enumerate}[label=\textnormal{(\roman*)}]
		\item There is a set of nonnegative multipliers $\{\nu_{i,j}\}_{i,j\in I_N^*}$ satisfying \eqref{eq:d1-tau}, such that the inequality \eqref{eq:G2-rate-modified} holds for all [families $\{y_i\}_{i\in I_N^*}$ in $\mathbb{R}^d$, $\{g_i\}_{i\in I_N^*}$ in $\mathbb{R}^d$, and $\{f_i\}_{i\in I_N^*}$ in $\mathbb{R}$ satisfying $g_*=0$ and \eqref{eq:fsfom-dual}]. 
		\item There is a set of nonnegative multipliers $\{\lambda_{i,j}\}_{i,j\in I_N^*}$ such that the inequality \eqref{eq:F2-rate} holds for all [families $\{x_i\}_{i\in I_N^*}$ in $\mathbb{R}^d$, $\{g_i\}_{i\in I_N^*}$ in $\mathbb{R}^d$, and $\{f_i\}_{i\in I_N^*}$ in $\mathbb{R}$ satisfying $g_*=0$ and \eqref{eq:fsfom}].
	\end{enumerate}
\end{theorem}
\section{Proof of the ``mirror'' of the main theorem}
\label{sec:pep}

In this section, we reprove Theorem~\ref{thm:gd}. Precisely, we show that there is a set of nonnegative multipliers $\{\nu_{i,j}\}_{i,j\in I_N^*}$ satisfying \eqref{eq:d1-tau}, such that \eqref{eq:G2-rate-modified} holds with 
\begin{equation}
	\label{eq:tau_value}
	\begin{aligned}
		\tau & =L\max\left\{ \frac{1}{1+\gamma L\sum_{j=1}^{2N}(1-\gamma\mu)^{-j}},(1-\gamma L)^{2N}\right\} \\
		& =L\max\left\{ \frac{1}{1+\gamma L\cdot E_{N}(1-\gamma\mu)},\frac{1}{1+\gamma L\cdot E_{N}(1-\gamma L)}\right\} 
	\end{aligned}
\end{equation}
for all [families $\{y_i\}_{i\in I_N^*}$ in $\mathbb{R}^d$, $\{g_i\}_{i\in I_N^*}$ in $\mathbb{R}^d$, and $\{f_i\}_{i\in I_N^*}$ in $\mathbb{R}$ satisfying $g_*=0$ and \eqref{eq:fsfom-dual} with $\bfH=(\gamma L)\bfI_N$]. Throughout this section, we assume that the value of $\tau$ is always given by \eqref{eq:tau_value}.

Note that if \eqref{eq:d1-tau} holds, then \eqref{eq:simple-g-condition} also holds. In this case, the inequality \eqref{eq:G2-rate-modified} can be equivalently written as
\begin{align*}
	& \left(\frac{1}{2\tau}-\frac{1}{2L}\right)\left\Vert g_{N}\right\Vert ^{2}-\frac{\mu L}{2(L-\mu)}\left\Vert y_{N}-y_{*}-\frac{1}{L}g_{N}\right\Vert ^{2}-f_{0}+f_{N}\\
	& \qquad\leq\sum_{i,j\in I_{N}}\nu_{i,j}\bigg[f_{j}-f_{i}+\left\langle g_{j},y_{i}-y_{j}\right\rangle +\frac{1}{2L}\left\Vert g_{i}-g_{j}\right\Vert ^{2}\\
	& \qquad\qquad\qquad\qquad\quad+\frac{\mu L}{2(L-\mu)}\left\Vert y_{i}-y_{j}-\frac{1}{L}\left(g_{i}-g_{j}\right)\right\Vert ^{2}\bigg].
\end{align*}
Thus, when $\mu\geq0$, a sufficient condition for \eqref{eq:G2-rate-modified} to hold is that the inequality
\begin{equation}
	\label{eq:sum_interpolation}
	\begin{aligned}
		& \left(\frac{1}{2\tau}-\frac{1}{2L}\right)\left\Vert g_{N}\right\Vert ^{2}-f_{0}+f_{N}\\
		& \qquad\leq\sum_{i,j\in I_{N}}\nu_{i,j}\bigg[f_{j}-f_{i}+\left\langle g_{j},y_{i}-y_{j}\right\rangle +\frac{1}{2L}\left\Vert g_{i}-g_{j}\right\Vert ^{2}\\
		& \qquad\qquad\qquad\qquad\quad+\frac{\mu L}{2(L-\mu)}\left\Vert y_{i}-y_{j}-\frac{1}{L}\left(g_{i}-g_{j}\right)\right\Vert ^{2}\bigg]
	\end{aligned}
\end{equation}
holds. Note that after substituting the value \eqref{eq:tau_value} of $\tau$, we have
\begin{equation}
	\label{eq:tauminL}
	\frac{1}{2\tau}-\frac{1}{2L}=\frac{\gamma}{2}\min\left\{ E_{N}(1-\gamma\mu),E_{N}(1-\gamma L)\right\} .
\end{equation}
To remind the goal, we want to show that there is a set of nonnegative multipliers $\{\nu_{i,j}\}_{i,j\in I_N}$ satisfying 
\begin{equation}
	\label{eq:d1-relax}
	\nu_{i,j}=0\quad\textup{for all }(i,j)\notin\left\{ (k,k+1),(k+1,k),(N,k)\right\} _{0\leq k\leq N-1},
\end{equation}
such that the inequality \eqref{eq:sum_interpolation} holds for all [families $\{y_i\}_{i\in I_N}$ in $\mathbb{R}^d$, $\{g_i\}_{i\in I_N}$ in $\mathbb{R}^d$, and $\{f_i\}_{i\in I_N}$ in $\mathbb{R}$ satisfying \eqref{eq:fsfom-dual} with $\bfH=(\gamma L)\bfI_N$].

As in Theorem~\ref{thm:G2-modified}, we can derive two conditions that together imply this goal. The first condition, which plays the role of \eqref{item:G2-1}, will be 
\begin{equation}
	\label{eq:linear_constraint}
	\sum_{i\in I_N}\nu_{i,k}-\sum_{j\in I_N}\nu_{k,j}=\begin{cases}
		1 & \textup{if }k=N,\\
		-1 & \textup{if }k=0,\\
		0 & \textup{otherwise}.
	\end{cases}    
\end{equation}
The second condition can be formulated similarly to \eqref{item:G2-2}, but we formulate it more explicitly in terms of the positive semidefiniteness of a certain matrix in the next subsection.

\subsection{Semidefinite reformulation of PEP analysis}
\label{sec:detail_pep}
Note that the inequality \eqref{eq:sum_interpolation} holds for all [families $\{y_i\}_{i\in I_N}$ in $\mathbb{R}^d$, $\{g_i\}_{i\in I_N}$ in $\mathbb{R}^d$, and $\{f_i\}_{i\in I_N}$ in $\mathbb{R}$ satisfying \eqref{eq:fsfom-dual}] if and only if the inequality
\begin{equation}
	\label{eq:sum_interpolation-2}
	\begin{aligned}
		& -\frac{1}{2L}\left\Vert g_{0}\right\Vert ^{2}+\frac{1}{2\tau}\left\Vert g_{N}\right\Vert ^{2}-f_{0}^+ +f_{N}^+\\
		& \qquad\leq\sum_{i,j\in I_{N}}\nu_{i,j}\left[f_{j}^{+}-f_{i}^{+}+\left\langle g_{j},y_{i}^{+}-y_{j}^{+}\right\rangle +\frac{\mu L}{2(L-\mu)}\left\Vert y_{i}^{+}-y_{j}^{+}\right\Vert ^{2}\right]
	\end{aligned}
\end{equation}
holds for all [families $\{y_0\}\cup\{y_i^+\}_{i\in I_N}$ in $\mathbb{R}^d$, $\{g_i\}_{i\in I_N}$ in $\mathbb{R}^d$, and $\{f_i^+\}_{i\in I_N}$ in $\mathbb{R}$ satisfying \eqref{eq:fsfom-dual-tilde}] (see Footnote~\ref{footnote:F}). Define two $(N+1)\times(N+1)$ matrices $\bfA$ and $\bfB$ as
\[
\bfA =\left[\begin{array}{ccccccc}
	0 & b_{0} & c_{0} & c_{0} & \cdots & c_{0} & c_{0}\\
	0 & a_{1} & b_{1} & c_{1} & \cdots & c_{1} & c_{1}\\
	0 & 0 & a_{2} & b_{2} & \cdots & c_{2} & c_{2}\\
	0 & 0 & 0 & a_{3} & \cdots & c_{3} & c_{3}\\
	\vdots & \vdots & \vdots & \vdots & \ddots & \vdots & \vdots\\
	0 & 0 & 0 & 0 & \cdots & a_{N-1} & b_{N-1}\\
	0 & 0 & 0 & 0 & \cdots & 0 & a_{N}
\end{array}\right],
\]
where
\begin{align*}
	a_{k} & =\nu_{k-1,k}\quad\textup{for }k=1,\ldots,N,\\
	b_{k} & =\begin{cases}
		-\nu_{k+1,k}-\nu_{N,k} & \textup{if }k\leq N-2,\\
		-\nu_{N,N-1} & \textup{if }k=N-1,
	\end{cases}\\
	c_{k} & =-\nu_{N,k}\quad\textup{for }k=0,\ldots,N-2,
\end{align*}
and
\[
\bfB =\left[\begin{array}{ccccccc}
	0 & 0 & 0 & 0 & \cdots & 0 & 0\\
	0 & d_{1} & e_{1} & e_{1} & \cdots & e_{1} & e_{1}\\
	0 & e_{1} & d_{2} & e_{2} & \cdots & e_{2} & e_{2}\\
	0 & e_{1} & e_{2} & d_{3} & \cdots & e_{3} & e_{3}\\
	\vdots & \vdots & \vdots & \vdots & \ddots & \vdots & \vdots\\
	0 & e_{1} & e_{2} & e_{3} & \cdots & d_{N-1} & e_{N-1}\\
	0 & e_{1} & e_{2} & e_{3} & \cdots & e_{N-1} & d_{N}
\end{array}\right],
\]
where
\begin{align*}
	d_{k} & =\begin{cases}
		\nu_{k-1,k}+\nu_{k,k-1}+\sum_{j=0}^{k-1}\nu_{N,j} & \textup{if }k\leq N-1,\\
		\nu_{N-1,N}+\sum_{j=0}^{N-1}\nu_{N,j} & \textup{if }k=N,
	\end{cases}\\
	e_{k} & =\sum_{j=0}^{k-1}\nu_{N,j}\quad\textup{for }k=0,\ldots,N-1.
\end{align*}
Then, we have
\begin{align*}
	\sum_{i,j\in I_N}\nu_{ij}\left\langle g_{j},y_{i}^{+}-y_{j}^{+}\right\rangle  &  =-\left\langle \bfg,\bfA\bfy\right\rangle , \\
	\sum_{i,j\in I_N}\nu_{ij}\left\Vert y_{i}^{+}-y_{j}^{+}\right\Vert ^{2} & =\left\langle \bfB\bfy,\bfy\right\rangle ,
\end{align*}
where $\bfg=[g_{0},\ldots,g_{N}]^{\rmT}$ and $\bfy=[y_{0}^{+}-y_{0},y_{1}^{+}-y_{0}^{+},\ldots,y_{N}^{+}-y_{N-1}^{+}]^{\rmT}$. Thus, when \eqref{eq:linear_constraint} holds, using $\bfg=-L\tilde{\bfH}^{-\rmA}\bfy$, the inequality \eqref{eq:sum_interpolation-2} can be equivalently written as
\begin{equation}
	\label{eq:appG2-2-A}
	0\leq S:=\frac{1}{2L}\left\Vert g_{0}\right\Vert ^{2}-\frac{1}{2\tau}\left\Vert g_{N}\right\Vert ^{2}+\left\langle \left(L\bfA^{\rmT}\tilde{\bfH}^{-\rmA}+\frac{\mu L}{2(L-\mu)}\bfB\right)\bfy,\bfy\right\rangle .
\end{equation}
Denote the first row and the last row of $\tilde{\bfH}^{-\rmA}$ by $\bfh_{0}$ and $\bfh_N$, both written as column vectors. Then, we have $S=L\left\langle \bfS\bfy,\bfy\right\rangle $, where
\begin{equation}
	\label{eq:pep-matrix}
	\bfS=\bfA^{\rmT}\tilde{\bfH}^{-\rmA}+\frac{\kappa}{2(1-\kappa)}\bfB+\frac{1}{2}\bfh_{0}\bfh_{0}^{\rmT}-\frac{L}{2\tau}\bfh_{N}\bfh_{N}^{\rmT}.
\end{equation}
Thus, we obtain a condition which plays the role of \eqref{item:G2-2}, namely the positive semidefiniteness of $\frac{1}{2}(\bfS+\bfS^{\rmT})$. We summarize the result in the following theorem.
\begin{theorem}
	\label{thm:pep}
	Given a set of nonnegative multipliers $\{\nu_{i,j}\}_{i,j\in I_N}$ satisfying \eqref{eq:d1-relax}, the inequality \eqref{eq:sum_interpolation} holds for all [families $\{y_i\}_{i\in I_N}$ in $\mathbb{R}^d$, $\{g_i\}_{i\in I_N}$ in $\mathbb{R}^d$, and $\{f_i\}_{i\in I_N}$ in $\mathbb{R}$ satisfying \eqref{eq:fsfom-dual}] if the following two conditions are satisfied:
	\begin{enumerate}[label=\textnormal{(G\arabic*$'$)},ref=G\arabic*$'$]
		\item \label{item:G3-1} The condition \eqref{eq:linear_constraint} holds.
		\item \label{item:G3-2} The matrix $\frac{1}{2}(\bfS+\bfS^{\rmT})$ is positive semidefinite, where $\bfS$ is defined in \eqref{eq:pep-matrix}.
	\end{enumerate}
\end{theorem}

\subsection{Proof of Theorem~\ref{thm:gd}}
\label{sec:proof_special}

Given $N$, $\mu$, and $L$ (possibly with $\mu<0$), the \emph{optimal stepsize} $\gamma^*(N,\mu,L)$ is uniquely defined as the value of $\gamma$ that optimizes the convergence rate \eqref{eq:gd-rate2}. It is known that $\gamma^*(N, \mu, L)$ is the unique value of $\gamma$ in $(0,2/L)$ that makes the two arguments of the $\min$ function in \eqref{eq:tauminL} equal. It is also known that $E_N(\eta) < E_N(\rho)$ for $\gamma \in(0, \gamma^*(N, \mu, L))$ and $E_N(\eta) > E_N(\rho)$ for $\gamma \in( \gamma^*(N, \mu, L),2/L)$. See \cite[\S2]{rotaru2024exact} for the mentioned facts. The following lemma shows the existence of a desired set of multipliers $\{\nu_{i,j}\}_{i,j\in I_N}$ for the case $\gamma=\gamma^*(N,\mu,L)$. 
\begin{lemma}
	\label{lem:pep-feasible}
	When $\mu\in(-\infty,L)$ and $\gamma=\gamma^*(N,\mu,L)$, there is a set of multipliers $\{\nu_{i,j}\}_{i,j\in I_N}$ satisfying \eqref{eq:d1-relax} such that the inequality \eqref{eq:sum_interpolation}, with $\tau$ set as in \eqref{eq:tau_value}, holds for all [families $\{y_i\}_{i\in I_N}$ in $\mathbb{R}^d$, $\{g_i\}_{i\in I_N}$ in $\mathbb{R}^d$, and $\{f_i\}_{i\in I_N}$ in $\mathbb{R}$ satisfying \eqref{eq:fsfom-dual} with $\bfH=(\gamma L)\bfI_N$].
\end{lemma}
In Appendix~\ref{app:lem:pep-feasible}, we prove Lemma~\ref{lem:pep-feasible} by explicitly setting the multipliers $\nu_{i,j}$ and then verifying the two conditions in Theorem~\ref{thm:pep}. It is nice that once Lemma~\ref{lem:pep-feasible} is proved, its generalization (Lemma~\ref{lem:pep-feasible-general}) follows effortlessly. The key idea is to consider $\gamma$ and $N$ as fixed, while allowing $\mu$ and $L$ to vary. 
\begin{lemma}
	\label{lem:pep-feasible-general}
	Lemma~\ref{lem:pep-feasible} holds even when the assumption $\gamma=\gamma^*(N,\mu,L)$ is replaced by $\gamma\in (0,2/L)$.
\end{lemma}
\begin{proof}
	Consider $\gamma$ and $N$ as fixed. Denote the set of multipliers given by Lemma~\ref{lem:pep-feasible} as $\{\nu_{i,j}\}_{i,j\in I_N}^{\mu,L}$, showing their dependency on $\mu$ and $L$.
	
	{\it Case}  1. $\gamma = \gamma^*(N,\mu,L)$. This is proved by Lemma~\ref{lem:pep-feasible}.
	
	{\it Case}  2. $\gamma < \gamma^*(N,\mu,L)$. In this case, we have $E_N(1-\gamma\mu)<E_N(1-\gamma L)$, meaning that $\mu$ determines the convergence rate. By the intermediate value theorem, we can choose $L'\in(L,2/\gamma)$ such that $E_N(1-\gamma\mu)=E_N(1-\gamma L')$ (see Table~\ref{table:1}), or equivalently $\gamma=\gamma^*(N,\mu,L')$. Consider the set of multipliers $\{\nu_{i,j}\}_{i,j\in I_N}^{\mu,L'}$ obtained from Lemma~\ref{lem:pep-feasible}. Then, for all [families $\{y_i\}_{i\in I_N}$ in $\mathbb{R}^d$, $\{g_i\}_{i\in I_N}$ in $\mathbb{R}^d$, and $\{f_i\}_{i\in I_N}$ in $\mathbb{R}$ satisfying \eqref{eq:fsfom-dual} with $\bfH=(\gamma L)\bfI_N$], the inequality \eqref{eq:sum_interpolation}, with $L$ replaced by $L'$, holds. By applying Proposition~\ref{prop:ineq-stronger}\eqref{item:L-compared} to each summand in \eqref{eq:sum_interpolation}, we can conclude that the inequality \eqref{eq:sum_interpolation} itself holds.
	
	{\it Case}  3. $\gamma > \gamma^*(N,\mu,L)$. In this case, we have $E_N(1-\gamma\mu)>E_N(1-\gamma L)$, meaning that $L$ determines the convergence rate. By the intermediate value theorem, we can choose $\mu' \in (-\infty,\mu)$ such that $E_N(1-\gamma\mu')=E_N(1-\gamma L)$ (see Table~\ref{table:1}), or equivalently $\gamma=\gamma^*(N,\mu',L)$. Consider the set of multipliers $\{\nu_{i,j}\}_{i,j\in I_N}^{\mu',L}$ obtained from Lemma~\ref{lem:pep-feasible}. Then, for all [families $\{y_i\}_{i\in I_N}$ in $\mathbb{R}^d$, $\{g_i\}_{i\in I_N}$ in $\mathbb{R}^d$, and $\{f_i\}_{i\in I_N}$ in $\mathbb{R}$ satisfying \eqref{eq:fsfom-dual} with $\bfH=(\gamma L)\bfI_N$], the inequality \eqref{eq:sum_interpolation}, with $\mu$ replaced by $\mu'$, holds. By applying Proposition~\ref{prop:ineq-stronger}\eqref{item:mu-compared} to each summand in \eqref{eq:sum_interpolation}, we can conclude that the inequality \eqref{eq:sum_interpolation} itself holds.
\end{proof}

\begin{proof}[Proof of Theorem~\ref{thm:gd}]
	Put $y_i\leftarrow x_i$ (see Footnote~\ref{footnote:abuse}), $g_i\leftarrow\nabla f(x_i)$, and $f_i\leftarrow f(x_i)$ for all $i\in I_N^*$. Then, the inequality \eqref{eq:sum_interpolation} holds by Lemma~\ref{lem:pep-feasible-general}. As discussed at the beginning of this section, this implies that the inequality \eqref{eq:G2-rate-modified} holds. Since the right-hand side of this inequality is nonpositive by Proposition~\ref{prop:int-ineq}, we can conclude that the desired inequality \eqref{eq:gd-rate2} holds.
\end{proof}

\begin{remark}
	Throughout this section, the condition $\mu\geq0$ was used only when showing that the inequality \eqref{eq:sum_interpolation} implies the inequality \eqref{eq:G2-rate-modified}. However, we can circumvent this step by using Proposition~\ref{prop:trivial}, which does not require $\mu\geq0$. Thus, the proof of Theorem~\ref{thm:gd} remains valid even when $\mu < 0$.
\end{remark}

\section{Proof of the main theorem}
\label{sec:gd}

We assume $\mu\geq0$.

\begin{proof}[Proof of Theorem~\ref{conj:gd}]
	For given $N$, $\mu$, $L$, and $\gamma$, consider the set of multipliers $\{\nu_{i,j}\}_{i,j\in I_N}$ given by Lemma~\ref{lem:pep-feasible-general}. Incorporate the index $*$ by setting the multipliers $\nu_{i,*}$ and $\nu_{*,j}$ as in \eqref{eq:simple-g-condition}. Then, as discussed in the beginning of Section~\ref{sec:pep}, the set of multipliers $\{\nu_{i,j}\}_{i,j\in I_N^*}$ satisfies \eqref{eq:d1-tau}, and that the inequality \eqref{eq:G2-rate-modified}, with $\tau$ set as in \eqref{eq:tau_value}, holds for all [families $\{y_i\}_{i\in I_N^*}$ in $\mathbb{R}^d$, $\{g_i\}_{i\in I_N^*}$ in $\mathbb{R}^d$, and $\{f_i\}_{i\in I_N^*}$ in $\mathbb{R}$ satisfying $g_*=0$ and \eqref{eq:fsfom-dual} with $\bfH=(\gamma L)\bfI_N$]. By Theorem~\ref{thm:relaxed_general}, there is a set of nonnegative multipliers $\{\lambda_{i,j}\}_{i,j\in I_N^*}$ such that the inequality \eqref{eq:F2-rate}, with $\tau$ set as in \eqref{eq:tau_value}, holds for all [families $\{x_i\}_{i\in I_N^*}$ in $\mathbb{R}^d$, $\{g_i\}_{i\in I_N^*}$ in $\mathbb{R}^d$, and $\{f_i\}_{i\in I_N^*}$ in $\mathbb{R}$ satisfying $g_*=0$ and \eqref{eq:fsfom} with $\bfH=(\gamma L)\bfI_N$]. We now put $x_i\leftarrow x_i$ (see Footnote~\ref{footnote:abuse}), $g_i\leftarrow\nabla f(x_i)$, and $f_i\leftarrow f(x_i)$ for all $i\in I_N^*$. Since the right-hand side of this inequality is nonpositive by Proposition~\ref{prop:int-ineq}, we can conclude that the desired inequality \eqref{eq:gd-rate} holds.
\end{proof}
Since our proofs of Lemmas~\ref{lem:relaxed_general} and \ref{lem:pep-feasible-general} are constructive, we can explicitly compute the set of multipliers $\{\lambda_{i,j}\}_{i,j\in I_N^*}$ used in the proof of Theorem~\ref{conj:gd}.
\begin{example}
	For $N=5$, $\mu=0.1$, $L=1$, and $\gamma=\gamma^*(N,\mu,L)$, the set of multipliers $\{\lambda_{i,j}\}_{i,j\in I_N^*}$ used in the proof of Theorem~\ref{conj:gd} is calculated as
	\begin{equation}
		\label{eq:lambda_exam}
		[\lambda_{i,j}] = \begin{blockarray}{cccccccc}
			& * & 0&1&2&3&4&5 \\
			\begin{block}{c[ccccccc]}
				* & 0 & 0.0384 & 0.0621 & 0.1063 & 0.1873 & 0.3342 & 0.2718\\
				0 & 0 & 0 & 0.0182 & 0.0119 & 0.0060 & 0.0020 & 0.0003\\
				1 & 0 & 0 & 0 & 0.0472 & 0.0237 & 0.0080 & 0.0013\\
				2 & 0 & 0 & 0 & 0 & 0.1186 & 0.0401 & 0.0067\\
				3 & 0 & 0 & 0 & 0 & 0 & 0.2876 & 0.0479\\
				4 & 0 & 0 & 0 & 0 & 0 & 0 & 0.6719\\
				5 & 0 & 0 & 0 & 0 & 0 & 0 & 0 \\
			\end{block}
		\end{blockarray}.
	\end{equation}
	A MATLAB code for generating the multipliers $\lambda_{i,j}$ can be found in the file \texttt{gener\-ate\_lambda.m} of the GitHub repository (see Section~\ref{sec:our_strategy}). 
\end{example}

\begin{remark}
	\label{rmk:grimmer}
	For the set of multipliers $\{\lambda_{i,j}\}_{i,j\in I_N^*}$ we obtained, one can see that
	\begin{equation}
		\label{eq:grimmer-property}
		{\lambda_{i,j}}/{\lambda_{i',j}}={\lambda_{i,j'}}/{\lambda_{i',j'}}
	\end{equation}
	holds whenever $0\leq i,i' < j,j'\leq N$. This can be proved from the explicit formula for converting $\{\nu_{i,j}\}_{i,j\in I_N^*}$ into $\{\lambda_{i,j}\}_{i,j\in I_N^*}$ provided in Appendix~\ref{app:lem:relaxed_general}. A similar property is satisfied by the conjectured set of multipliers in \cite[Conj.~2]{grimmer2024strengthened}.\footnote{In \cite{grimmer2024strengthened}, it was conjectured that there is a set of multipliers $\{\lambda_{i,j}\}_{i,j\in I_N^*}$ satisfying \eqref{eq:grimmer-property} whenever $i,i' < j-1,j'-1$, which proves Theroem~\ref{conj:gd} under a PEP analysis which is similar to that presented in Section~\ref{sec:F2}. The structure of their $\{\lambda_{i,j}\}_{i,j \in I_N^*}$ is different from ours.}
\end{remark}

\appendix

\section{Proof of Lemma~\ref{lem:relaxed_general}}
\label{app:lem:relaxed_general}

Let $A_{k}=1+\sum_{j=0}^{k-1}\nu_{N,j}$ for $k=1,\ldots,N$ and $B_k=\nu_{k,k-1}$ for $k=1,\ldots,N-1$. For convenience, let $B_N=0$. Then, $\{A_k\}_{k=1}^N$ is a nondecreasing sequence in $[1,\infty)$ and $\{B_k\}_{k=1}^N$ is a sequence in $[0,\infty)$. Using \eqref{eq:relaxed-g} and \eqref{item:G2-1}, we can express all the multipliers $\nu_{i,j}$ as 
\[
\nu_{i,j}=\begin{cases}
	A_{j}+B_{j} & \textup{if }i=j-1,\ 1\leq j\leq N,\\
	A_{1}-1 & \textup{if }i=N,\ j=0,\\
	A_{j+1}-A_{j} & \textup{if }i=N,\ 1\leq j\leq N-1,\\
	B_{i} & \textup{if }j=i-1,\ 1\leq i\leq N-1,\\
	1 & \textup{if }i=N,\ j=*,\\
	0 & \textup{otherwise}.
\end{cases}
\]
Thus, \eqref{eq:ymapx2} without the last equality can be expressed as \eqref{eq:transform_matrix} with
\[
\bfM_{i,j}=\begin{cases}
	A_{i}+B_{i} & \textup{if }j=i,\\
	A_{i}-A_{i+1}-B_{i+1} & \textup{if }j=i+1,\\
	A_{i}-A_{i+1} & \textup{if } i+2\leq j\leq N.
\end{cases}.
\]
The last equality in \eqref{eq:ymapx2} can be expressed as \eqref{eq:last_simple}. The entries of the inverse matrix $\bfM^{-1}$ can be computed as
\begin{equation}
	\label{eq:general-Minv}
	\bfM_{i,j}^{-1}=\begin{cases}
		\frac{1}{A_{i}+B_{i}} & \textup{if }j=i,\\
		\frac{1}{A_{i}+B_{i}}-\frac{1}{A_{i}+B_{i}}\frac{A_{i}}{A_{i+1}+B_{i+1}} & \textup{if }j=i+1,\\
		\frac{1}{A_{i}+B_{i}}-\sum_{l=i+1}^{j}\big[\frac{1}{A_{i}+B_{i}}\frac{A_{i}}{A_{i+1}+B_{i+1}}\\
		\qquad\qquad\qquad\qquad\times\prod_{k=i+1}^{l-1}\frac{B_{k}}{A_{k+1}+B_{k+1}}\big] & \textup{if }i+2\leq j\leq N.
	\end{cases}
\end{equation}
A MATLAB code for verifying \eqref{eq:general-Minv} can be found in the file \texttt{inverse\_transform.m} of the GitHub repository (see Section~\ref{sec:our_strategy}). We set the multipliers $\lambda_{i,j}$ as
\[
\lambda_{i,j}=\begin{cases}
	\bfM_{N-j,N-i-1}^{-1}-\bfM_{N-j,N-i}^{-1} & \textup{if }0\leq i<j\leq N-1,\\
	\bfM_{N-j,N}^{-1} & \textup{if }i=*,\ 0\leq j\leq N-1,\\
	\bfM_{N-i,N-i}^{-1}-\sum_{k=1}^{N-i-1}\left(\bfM_{k,N-i-1}^{-1}-\bfM_{k,N-i}^{-1}\right) & \textup{if }0\leq i<j=N,\\
	1-\sum_{l=0}^{N-1}\Big[\bfM_{N-l,N-l}^{-1}\\
	\qquad\qquad\quad-\sum_{k=1}^{N-l-1}\left(\bfM_{k,N-l-1}^{-1}-\bfM_{k,N-l}^{-1}\right)\Big] & \textup{if }i=*,\ j=N,\\
	0 & \textup{otherwise}.
\end{cases}
\]
After substituting \eqref{eq:general-Minv}, we have
\[
\lambda_{i,j}=\begin{cases}
	\frac{1}{A_{N-j}+B_{N-j}}\frac{A_{N-j}}{A_{N-j+1}+B_{N-j+1}}\\
	\qquad\times\prod_{k=N-j+1}^{N-i-1}\frac{B_{k}}{A_{k+1}+B_{k+1}} & \textup{if }0\leq i<j\leq N-1,\\
	\frac{1}{A_{N-j}+B_{N-j}}-\sum_{l=0}^{j-1}\big[\frac{1}{A_{N-j}+B_{N-j}}\\
	\qquad\qquad\qquad\qquad\qquad\times\frac{A_{N-j}}{A_{N-j+1}+B_{N-j+1}}\\
	\qquad\qquad\qquad\qquad\qquad\times\prod_{k=N-j+1}^{N-l-1}\frac{B_{k}}{A_{k+1}+B_{k+1}}\big] & \textup{if }i=*,0\leq j\leq N-1,\\
	\frac{1}{A_{N-i}+B_{N-i}}\prod_{k=1}^{N-i-1}\frac{B_{k}}{A_{k}+B_{k}} & \textup{if }0\leq i<j=N,\\
	1-\sum_{l=0}^{N-1}\big[\frac{1}{A_{N-l}+B_{N-l}}\prod_{k=1}^{N-l-1}\frac{B_{k}}{A_{k}+B_{k}}\big] & \textup{if }i=*,j=N,\\
	0 & \textup{otherwise}.
\end{cases}
\]
Then, the condition \eqref{eq:xmapy2} is equivalent to \eqref{eq:last_simple} and \eqref{eq:transform-inverse} with \eqref{eq:general-Minv}. Thus, \eqref{eq:xmapy2} is equivalent to \eqref{eq:ymapx2}. It is straightforward to check \eqref{item:F2-1}.

It only remains to check $\lambda_{i,j}\geq 0$ for all $i,j\in I_N^*$. For the case where $0\leq i<j\leq N-1$ and the case where $i=*$ and $0\leq j\leq N-1$, this is clear. For the case where $i=*$ and $0\leq j\leq N-1$, we have
\begin{align*}
	\lambda_{i,j} & \geq\frac{1}{A_{N-j}+B_{N-j}}\\
	& \quad-\sum_{l=0}^{j-1}\left(\frac{1}{A_{N-j}+B_{N-j}}\frac{A_{N-j}}{A_{N-j}+B_{N-j+1}}\prod_{k=N-j+1}^{N-l-1}\frac{B_{k}}{A_{N-j}+B_{k+1}}\right)\\
	& =\frac{1}{A_{N-j}+B_{N-j}}\prod_{k=N-j+1}^{N}\frac{B_{k}}{A_{N-j}+B_{k}}\geq0,
\end{align*}
where for the first inequality we used the fact that $A_k \geq A_{N-j}$ whenever $k \geq N-j$. For the case where $i=*$ and $j=N$, we have
\[
\lambda_{i,j}\geq1-\sum_{l=0}^{N-1}\left[\frac{1}{1+B_{N-l}}\prod_{k=1}^{N-l-1}\frac{B_{k}}{1+B_{k}}\right]=\prod_{k=1}^{N}\frac{B_{k}}{1+B_{k}}\geq0,
\]
where for the first inequality we used the fact that $A_k\geq 1$ for all $k$. This completes the proof.

\section{Proof of Lemma~\ref{lem:pep-feasible}}
\label{app:lem:pep-feasible}

After choosing a set of multipliers $\{\nu_{i,j}\}_{i,j\in I_N}$, we verify the conditions \eqref{eq:d1-relax}, \eqref{eq:linear_constraint}, $\nu_{i,j} \geq 0$ for all $i,j\in I$, and that the matrix $\frac{1}{2}(\bfS+\bfS^{\rmT})$ is positive semidefinite. The following propositions will be used throughout the proof.
\begin{proposition}[{\cite[Props.~2.12, 4.6]{rotaru2024exact}}]
	\label{prop:tkpos}
	When $\gamma=\gamma^*(N,\mu,L)$, we have $\rho\in(-1,0)$, $\eta\in(-\rho,\infty)$, $T_k(\rho,\eta)\geq 0$ for $k=1,\ldots,N-1$, and $T_N(\rho,\eta)=0$.
\end{proposition}
\begin{proposition}
	\label{prop:convex}
	For any $\rho\in(-1,0)$, $\eta\in(0,\infty)$, and $N\in(0,\infty)$, the function
	\[
	\psi(t)=\begin{cases}
		\log\left(\frac{1+(1-\rho)(N+t)}{1+(1-\rho)(N-t)}\right) & \textup{if }\eta=1,\\
		\log\left(\frac{-(\eta-\rho)+(1-\rho)\eta^{-t-N}}{-(\eta-\rho)+(1-\rho)\eta^{t-N}}\right) & \textup{otherwise},
	\end{cases}
	\]
	is convex on $[0,N]$.
\end{proposition}

\begin{proof} 
	We present a proof for the case $\eta\neq 1$. The case $\eta=1$ can be proved by applying similar reasoning. We need to show that the derivative
	\[
	\frac{d\psi}{dt}=\frac{-\left(\log\eta\right)(1-\rho)\eta^{-t-N}}{-(\eta-\rho)+(1-\rho)\eta^{-t-N}}+\frac{-\left(\log\eta\right)(1-\rho)\eta^{t-N}}{-(\eta-\rho)+(1-\rho)\eta^{t-N}}
	\]
	is nondecreasing on $[0,N]$. Note that if a differentiable function $\varphi$ on $(-\infty,0]$ is convex, then the function $t\mapsto\varphi(-N+t)+\varphi(-N-t)$ on $[0,N]$ is nondecreasing.\footnote{We have $\frac{d}{dt}\{\varphi(-N+t)+\varphi(-N-t)\}=\frac{d\varphi}{dt}(-N+t)-\frac{d\varphi}{dt}(-N-t)\geq0$, where the inequality follows from the fact that $\frac{d\varphi}{dt}$ is nondecreasing.} Thus, it suffices to show that the function
	\[
	\varphi(u)=\frac{-(\log\eta)(1-\rho)\eta^{u}}{-(\eta-\rho)+(1-\rho)\eta^{u}}=-(\log\eta)\left(1+\frac{(\eta-\rho)}{-(\eta-\rho)+(1-\rho)\eta^{u}}\right)
	\]
	is convex on $(-\infty,0]$. Its second derivative can be calculated as
	\[
	\frac{d^{2}\varphi}{du^{2}}=\frac{\left(\log\eta\right)^{3}\eta^{u}(\eta-\rho)(1-\rho)\left(-(\eta-\rho)+(1-\rho)\eta^{u}\right)\left(-(\eta-\rho)-(1-\rho)\eta^{u}\right)}{\left(-(\eta-\rho)+(1-\rho)\eta^{u}\right)^{4}}.
	\]
	It is straightforward to see that $\frac{d^{2}\varphi}{du^{2}}(u)\geq0$ for all $u\in(-\infty,0)$, meaning that $\varphi$ is convex on $(-\infty,0]$. This completes the proof.
\end{proof}

\subsection{Choosing multipliers}
\label{app:choosing_multipliers}

In this subsection, we set multipliers $\nu_{i,j}$. Since $\gamma=\gamma^*(N,\mu,L)$, we have $\rho,\eta\neq 0$ by Proposition~\ref{prop:tkpos}. Thus, $E_k(\eta)$ and $E_k(\rho)$ are finite for $k=1,\ldots,N$. 

{\it Case}  1. $N=1$. We set the multipliers $\nu_{i,j}$ as
\[
\nu_{i,j}=\begin{cases}
	-\rho E_{1}(\rho) & \textup{if }i=1,\ j=0,\\
	1-\rho E_{1}(\rho) & \textup{if }i=0,\ j=1.
\end{cases}
\]

{\it Case}  2. $N=2$. Define two constants $\alpha_1$ and $\beta_1$ as
\[
\alpha_{1} =\frac{T_{1}(\rho,\eta)}{F_{1}(\eta)}, \qquad \beta_{1} =\frac{\eta-\rho}{\eta}E_{1}(\rho)-\frac{T_{1}(\rho,\eta)}{F_{1}(\eta)}.
\]
We set the multipliers $\nu_{i,j}$ as
\[
\nu_{i,j}=\begin{cases}
	-\frac{\eta\rho}{\eta-\rho}\beta_{1} & \textup{if }i=1,\ j=0,\\
	-\frac{\eta\rho}{\eta-\rho}\alpha_{1} & \textup{if }i=2,\ j=0,\\
	1-\frac{\eta\rho}{\eta-\rho}\left(\alpha_{1}+\beta_{1}\right) & \textup{if }i=0,\ j=1,\\
	-\rho E_{2}(\rho)+\frac{1+\rho}{\eta-\rho}T_{1}(\rho,\eta) & \textup{if }i=2,\ j=1,\\
	1-\frac{\eta\rho}{\eta-\rho}\alpha_{1}-\rho E_{2}(\rho)+\frac{1+\rho}{\eta-\rho}T_{1}(\rho,\eta) & \textup{if }i=1,\ j=2,\\
	0 & \textup{otherwise}.
\end{cases}
\]

{\it Case}  3. $N\geq 3$. Define two sequences $\{\alpha_k\}_{k=1}^{N-1}$ and $\{\beta_k\}_{k=1}^{N-1}$ as
\[
\alpha_{k}=\begin{cases}
	\frac{T_{1}(\rho,\eta)}{F_{N-1}(\eta)} & \textup{if }k=1,\\
	-\frac{1}{\rho}\frac{T_{1}(\rho,\eta)}{F_{N-1}(\eta)}+\left(\frac{T_{2}(\rho,\eta)}{F_{N-2}(\eta)}-\frac{T_{1}(\rho,\eta)}{F_{N-1}(\eta)}\right) & \textup{if }k=2,\\
	-\frac{1}{\rho}\left(\frac{T_{k-1}(\rho,\eta)}{F_{N-k+1}(\eta)}-\frac{T_{k-2}(\rho,\eta)}{F_{N-k+2}(\eta)}\right)+\left(\frac{T_{k}(\rho,\eta)}{F_{N-k}(\eta)}-\frac{T_{k-1}(\rho,\eta)}{F_{N-k+1}(\eta)}\right) & \textup{if }3\leq k\leq N-1,
\end{cases}
\]
and
\[
\beta_{k} =\frac{\eta-\rho}{\eta}E_{k}(\rho)-\frac{T_{k}(\rho,\eta)}{F_{N-k}(\eta)}\quad\textup{for }k=1,\ldots,N-1.
\]
We set the multipliers $\nu_{i,j}$ as
\[
\nu_{i,j}=\begin{cases}
	-\frac{\eta\rho}{\eta-\rho}\beta_{j+1} & \textup{if }i=j+1,\ 0\leq j\leq N-2,\\
	-\frac{\eta\rho}{\eta-\rho}\alpha_{j+1} & \textup{if }i=N,\ 0\leq j\leq N-2,\\
	1-\frac{\eta\rho}{\eta-\rho}\left(\sum_{k=1}^{j}\alpha_{k}+\beta_{j}\right) & \textup{if }i=j-1,\ 1\leq j\leq N-1,\\
	\frac{\eta}{\eta-\rho}\left(\frac{T_{N-1}(\rho,\eta)}{F_{1}(\eta)}-\frac{T_{N-2}(\rho,\eta)}{F_{2}(\eta)}\right)\\
	\qquad-\rho E_{N}(\rho)+\frac{\eta\rho}{\eta-\rho}\frac{T_{N-1}(\rho,\eta)}{F_{1}(\eta)} & \textup{if }i=N,\ j=N-1,\\
	\frac{\eta}{\eta-\rho}\left(\frac{T_{N-1}(\rho,\eta)}{F_{1}(\eta)}-\frac{T_{N-2}(\rho,\eta)}{F_{2}(\eta)}\right)-\rho E_{N}(\rho)\\
	\qquad+\frac{\eta\rho}{\eta-\rho}\frac{T_{N-1}(\rho,\eta)}{F_{1}(\eta)}-\frac{\eta\rho}{\eta-\rho}\sum_{k=1}^{N-1}\alpha_{k}+1 & \textup{if }i=N-1,\ j=N,\\
	0 & \textup{otherwise}.
\end{cases}
\]

It is clear that the conditions \eqref{eq:d1-relax} and \eqref{eq:linear_constraint} hold for any case.

\subsection{Verifying nonnegativity of multipliers}
\label{app:verifying-nonneg}

For simplicity, we only provide a proof for the third case in the previous subsection. The second case is similar. The first case is trivial. We only need to show that $\alpha_k$, $\beta_k$, and $\nu_{N,N-1}$ are nonnegative, since every multiplier $\nu_{i,j}$ is a weighted sum of these plus a nonnegative constant.

\subsubsection{Verifying nonnegativity of $\alpha_k$}
\label{app:a1}

For $k=1,\ldots,N-1$, the nonnegativity of $\alpha_k$ follows from Proposition~\ref{prop:tkpos} and the following proposition.
\begin{proposition}
	\label{prop:for_alpha}
	The inequality
	\begin{equation}
		\label{eq:for_alpha}
		\frac{T_{k+1}(\rho,\eta)}{F_{N-k-1}(\eta)}-\frac{T_{k}(\rho,\eta)}{F_{N-k}(\eta)}\geq0
	\end{equation}
	holds for $k=1,\ldots,N-2$
\end{proposition}
\begin{proof}
	We consider two cases: $\eta\in(0,1]$ and $\eta\in(1,\infty)$.
	
	{\it Case}  1. $\eta\in(0,1]$. Define a sequence $\{\varphi_k\}_{k=1}^{N}$ by
	\begin{align*}
		\varphi_{k} & =\eta^{2k-2}\left(\eta^{-2k+1}+\eta^{-2k}-\rho^{-2k+1}-\rho^{-2k}\right)\\
		& =\left(\eta^{-1}+\eta^{-2}\right)-\frac{\eta^{2k-2}}{\rho^{2k-2}}\left(\rho^{-1}+\rho^{-2}\right).
	\end{align*}
	Then, $\varphi_k$ is nonincreasing by Proposition~\ref{prop:tkpos}, and we have
	\[
	T_{k}(\rho,\eta)=\varphi_{1}+\eta^{-2}\varphi_{2}+\cdots+\eta^{-2k+2}\varphi_{k}.
	\]
	Thus, for $k=1,\ldots,N-2$, we have
	\begin{align*}
		0 & =T_{N}(\rho,\eta)\\
		& =\varphi_{1}+\eta^{-2}\varphi_{2}+\cdots+\eta^{-2N+2}\varphi_{N}\\
		& \leq\varphi_{1}+\eta^{-2}\varphi_{2}+\cdots+\eta^{-2k+2}\varphi_{k}\\
		& \quad+\left(\eta^{-2k}+\eta^{-2k-2}\cdots+\eta^{-2N+2}\right)\varphi_{k+1}\\
		& =T_{k}(\rho,\eta)+\left(\eta^{-2k}+\eta^{-2k-2}+\cdots+\eta^{-2N+2}\right)\varphi_{k+1}.
	\end{align*}
	Using $T_{k+1}(\rho,\eta)=T_{k}(\rho,\eta)+\eta^{-2k}\varphi_{k+1}$, we have
	\begin{align*}
		0 & \leq T_{k}(\rho,\eta)+\left(\eta^{-2k}+\eta^{-2k-2}\cdots+\eta^{-2N+2}\right)\eta^{2k}\left(T_{k+1}(\rho,\eta)-T_{k}(\rho,\eta)\right)\\
		& =\left(1+\eta^{-2}+\cdots+\eta^{-2N+2k+2}\right)T_{k+1}(\rho,\eta)-\left(\eta^{-2}+\cdots+\eta^{-2N+2k+2}\right)T_{k}(\rho,\eta).
	\end{align*}
	Thus, we have
	\[
	\frac{T_{k+1}(\rho,\eta)}{T_{k}(\rho,\eta)}\geq \frac{\eta^{-2}+\cdots+\eta^{-2N+2k+2}}{1+\eta^{-2}+\cdots+\eta^{-2N+2k+2}}.
	\]
	To prove \eqref{eq:for_alpha}, we only need to show
	\[
	\frac{\eta^{-2}+\cdots+\eta^{-2N+2k+2}}{1+\eta^{-2}+\cdots+\eta^{-2N+2k+2}}\geq\frac{\eta+\eta^{2}+\cdots+\eta^{N-k-1}}{\eta+\eta^{2}+\cdots+\eta^{N-k}}.
	\]
	One can verify that this inequality holds for all $\eta\in(0,1]$.
	
	{\it Case} 2. $\eta\in(1,\infty)$. Define a sequence $\{\varphi_k\}_{k=1}^{N}$ by
	\[
	\varphi_k=\eta^{-2k+1}+\eta^{-2k}-\rho^{-2k+1}-\rho^{-2k}.
	\]
	Then, $\varphi_k$ is nonincreasing by Proposition~\ref{prop:tkpos}, and we have
	\[
	T_{k}(\rho,\eta)=\varphi_{1}+\varphi_{2}+\cdots+\varphi_{k}.
	\]
	Following a similar argument as in the first case, we arrive at
	\[
	\frac{T_{k+1}(\rho,\eta)}{T_{k}(\rho,\eta)}\geq\frac{N-k-1}{N-k}.
	\]
	Thus, we only need to show
	\[
	\frac{N-k-1}{N-k}\geq\frac{\eta+\eta^{2}+\cdots+\eta^{N-k-1}}{\eta+\eta^{2}+\cdots+\eta^{N-k}}.
	\]
	One can verify that this inequality holds for all $\eta\in(1,\infty)$.
\end{proof}

\subsubsection{Verifying nonnegativity of $\beta_k$}
To show $\beta_k\geq0$ for $k=1,\ldots,N-1$, it suffices to prove the following proposition.
\begin{proposition}
	\label{prop:for_beta}
	The inequality
	\begin{equation}
		\label{eq:for_beta}
		\frac{\eta-\rho}{\eta}E_k(\rho)-\frac{T_k(\rho,\eta)}{F_{N-k}(\eta)}\geq 0
	\end{equation}
	holds for $k=1,\ldots,N-1$.
\end{proposition}
\begin{proof}
	We present a proof for the case where $\eta\neq 1$. The case where $\eta=1$ can be handled using the same argument with the expressions $E_k(\eta)=2k$ and $F_{N-k}(\eta)=N-k$. Substituting $T_k(\rho,\eta)=E_k(\eta)-E_k(\rho)$ and dividing both sides by $\frac{\eta-\rho}{\eta}+\frac{1}{F_{N-k}(\eta)}$, the inequality \eqref{eq:for_beta} can be equivalently written as
	\begin{equation}
		\label{eq:for_beta2-before}
		E_{k}(\rho)-\frac{\eta}{(\eta-\rho)F_{N-k}(\eta)+\eta}E_{k}(\eta)\geq0.
	\end{equation}
	Using the expressions $E_{k}(x)=\frac{x^{-2k}-1}{1-x}$ and $F_{k}(x)=\frac{x(1-x^{k})}{1-x}$, this inequality can be equivalently written as
	\begin{equation}
		\label{eq:for_beta2}
		\frac{(-\rho)^{-2k}-1}{1-\rho}\geq\frac{\eta^{-2k}-1}{(\eta-\rho)\left(1-\eta^{N-k}\right)+(1-\eta)}.
	\end{equation}
	Multiplying both sides by $1-\rho$ and adding $1$ to each side, we obtain
	\begin{align*}
		(-\rho)^{-2k} & \geq\frac{(1-\rho)\left(\eta^{-2k}-1\right)}{(\eta-\rho)\left(1-\eta^{N-k}\right)+(1-\eta)}+1\\
		& =\frac{-(\eta-\rho)+(1-\rho)\eta^{-k-N}}{-(\eta-\rho)+(1-\rho)\eta^{k-N}}.
	\end{align*}
	By taking logarithms, we obtain
	\[
	-2k\log(-\rho)\geq\log\left(\frac{-(\eta-\rho)+(1-\rho)\eta^{-k-N}}{-(\eta-\rho)+(1-\rho)\eta^{k-N}}\right).
	\]
	We consider $k$ as a real variable, and view both sides as functions of $k$. This inequality becomes an equality at $k=0$ and $k=N$, as is clear from its equivalent forms \eqref{eq:for_beta2-before} and \eqref{eq:for_beta2}. By Proposition~\ref{prop:convex}, the right-hand side is convex on $[0,N]$. Therefore, the given inequality is valid for all $k \in [0,N]$.
\end{proof}

\subsubsection{Verifying nonnegativity of $\nu_{N,N-1}$}

We rewrite $\nu_{N,N-1}$ as
\begin{align*}
	\nu_{N,N-1}& =-\rho E_{N}(\rho)+\frac{\rho}{\eta-\rho}T_{N-1}(\rho,\eta)+\frac{\eta}{\eta-\rho}\left(\frac{T_{N-1}(\rho,\eta)}{\eta}-\frac{T_{N-2}(\rho,\eta)}{\eta+\eta^{2}}\right)\\
	& =-\frac{\rho\eta}{\eta-\rho}\left(\frac{\eta-\rho}{\eta}E_{N-1}(\rho)-\frac{T_{N-1}(\rho,\eta)}{F_{1}(\eta)}\right)\\
	& \quad +\frac{\eta}{\eta-\rho}\left(\frac{T_{N-1}(\rho,\eta)}{F_{1}(\eta)}-\frac{T_{N-2}(\rho,\eta)}{F_{2}(\eta)}\right)-\rho\left(\rho^{-2N+1}+\rho^{-2N}\right). \\
\end{align*}
The nonnegativity of $\nu_{N,N-1}$ now follows from Propositions~\ref{prop:tkpos}, \ref{prop:for_alpha}, and \ref{prop:for_beta}.

\subsection{Verifying positive semidefiniteness of $\frac{1}{2}(\bfS+\bfS^{\rmT})$}
\label{app:a2}

From $\bfH=(\gamma L)\bfI_N=(1-\rho)\bfI_{N}$, we obtain
\[
\tilde{\bfH}^{-\rmA}=\left[\begin{array}{ccccc}
	1\\
	\rho & 1\\
	\rho^{2} & \rho & 1\\
	\vdots & \vdots & \vdots & \ddots\\
	\rho^{N} & \rho^{N-1} & \rho^{N-2} & \cdots & 1
\end{array}\right].
\]
Since $E_N(\eta)=E_N(\rho)$, we have $\tau=L \rho^{2N}$. We substitute all expressions into \eqref{eq:pep-matrix}, and then find a nice expression for the matrix $\frac{1}{2}(\bfS+\bfS^{\rmT})$. This task is summarized in the following proposition.
\begin{proposition}
	\label{prop:PEP-simple}
	With $\bfH=(1-\rho)\bfI_{N}$, $\tau=L\rho^{2N}$, and the multipliers $\nu_{i,j}$ chosen in Appendix~\ref{app:choosing_multipliers}, the matrix $\frac{1}{2}(\bfS+\bfS^{\rmT})$ in Section~\ref{sec:detail_pep} can be expressed as
	\begin{equation}
		\label{eq:PEP-matrix-simple}
		\frac{1}{2}(\bfS+\bfS^{\rmT})=\begin{cases}
			\frac{\eta^{2}(1-\rho)}{2(\eta-\rho)^{2}}\left[\begin{smallmatrix}
				0 & 0\\
				0 & T_{1}(\rho,\eta)
			\end{smallmatrix}\right] & \textup{if }N=1,\\
			\frac{\eta^{2}(1-\rho)}{2(\eta-\rho)^{2}}\sum_{k=1}^{N}\delta_{k}\bfv_{k}\bfv_{k}^{\rmT} & \textup{if }N\geq2,
		\end{cases}
	\end{equation}
	where the sequence $\{\delta_k\}_{k=1}^{N}$ is defined by
	\[
	\delta_{k}=\begin{cases}
		T_{1}(\rho,\eta) & \textup{if }k=1,\\
		T_{k}(\rho,\eta)-\frac{F_{N-k}(\eta)^{2}}{F_{N-k+1}(\eta)^{2}}T_{k-1}(\rho,\eta) & \textup{if }2\leq k\leq N,
	\end{cases}
	\]
	and the vectors $\bfv_1,\ldots,\bfv_N\in\mathbb{R}^{N+1}$ are defined by
	\begin{align*}
		\bfv_{1} & =\left[0,1,-1/F_{N-1}(\eta),\ldots,-1/F_{N-1}(\eta)\right]^{\rmT},\\
		\bfv_{2} & =\left[0,0,1,-1/F_{N-2}(\eta),\ldots,-1/F_{N-2}(\eta)\right]^{\rmT},\\
		& \vdotswithin{=}\\
		\bfv_{N-1} & =\left[0,\ldots,0,1,-1/F_{1}(\eta)\right]^{\rmT},\\
		\bfv_{N} & =\left[0,\ldots,0,1\right]^{\rmT}.
	\end{align*}
\end{proposition}
A MATLAB code for verifying Proposition~\ref{prop:PEP-simple} can be found in the file \texttt{pep\_mat\-rix.m} of the GitHub repository (see Section~\ref{sec:our_strategy}). Due to this nice expression, proving the positive semidefiniteness of $\frac{1}{2}(\bfS+\bfS^{\rmT})$ reduces to proving $\delta_k \geq 0$ for $k=1,\ldots,N$. This follows from $\frac{F_{N-k}(\eta)}{F_{N-k+1}(\eta)}\in[0,1]$, Propositions~\ref{prop:tkpos}, and \ref{prop:for_alpha}.

\bibliographystyle{amsplain}
\bibliography{gdexact}

\providecommand{\bysame}{\leavevmode\hbox to3em{\hrulefill}\thinspace}
\providecommand{\MR}{\relax\ifhmode\unskip\space\fi MR }
\providecommand{\MRhref}[2]{%
  \href{http://www.ams.org/mathscinet-getitem?mr=#1}{#2}
}
\providecommand{\href}[2]{#2}
\begin{thebibliography}{10}

\bibitem{bubeck2015convex}
Sébastien Bubeck, \emph{Convex optimization: Algorithms and complexity},
  Foundations and Trends® in Machine Learning \textbf{8} (2015), no.~3-4,
  231--357.

\bibitem{cauchy1847methode}
Augustin Cauchy et~al., \emph{M{\'e}thode g{\'e}n{\'e}rale pour la
  r{\'e}solution des systemes d’{\'e}quations simultan{\'e}es}, C. R. Acad.
  Sci. Paris \textbf{25} (1847), no.~1847, 536--538.

\bibitem{curry1944method}
Haskell~B. Curry, \emph{The method of steepest descent for non-linear
  minimization problems}, Quart. Appl. Math. \textbf{2} (1944), 258--261.
  \MR{10667}

\bibitem{drori2014performance}
Yoel Drori and Marc Teboulle, \emph{Performance of first-order methods for
  smooth convex minimization: a novel approach}, Math. Program. \textbf{145}
  (2014), no.~1-2, 451--482. \MR{3207695}

\bibitem{daspremont2021acceleration}
Alexandre d’Aspremont, Damien Scieur, and Adrien Taylor, \emph{Acceleration
  methods}, Foundations and Trends® in Optimization \textbf{5} (2021),
  no.~1-2, 1--245.

\bibitem{golyshev2007fuchsian}
Vasily Golyshev and Jan Stienstra, \emph{Fuchsian equations of type {DN}},
  Commun. Number Theory Phys. \textbf{1} (2007), no.~2, 323--346. \MR{2346574}

\bibitem{grimmer2024strengthened}
Benjamin Grimmer, Kevin Shu, and Alex~L Wang, \emph{A strengthened conjecture
  on the minimax optimal constant stepsize for gradient descent}, arXiv
  preprint arXiv:2407.11739 (2024).

\bibitem{hadamard1908memoire}
Jacques Hadamard, \emph{M{\'e}moire sur le probl{\`e}me d'analyse relatif {\`a}
  l'{\'e}quilibre des plaques {\'e}lastiques encastr{\'e}es}, vol.~33,
  Imprimerie nationale, 1908.

\bibitem{kim2021optimizing}
Donghwan Kim and Jeffrey~A. Fessler, \emph{Optimizing the efficiency of
  first-order methods for decreasing the gradient of smooth convex functions},
  J. Optim. Theory Appl. \textbf{188} (2021), no.~1, 192--219. \MR{4200948}

\bibitem{kim2024time}
Jaeyeon Kim, Asuman Ozdaglar, Chanwoo Park, and Ernest Ryu, \emph{Time-reversed
  dissipation induces duality between minimizing gradient norm and function
  value}, Advances in Neural Information Processing Systems \textbf{36} (2024).

\bibitem{kim2024convergence}
Jungbin Kim and Insoon Yang, \emph{Convergence analysis of ode models for
  accelerated first-order methods via positive semidefinite kernels}, Advances
  in Neural Information Processing Systems \textbf{36} (2024).

\bibitem{nesterov2018lectures}
Yurii Nesterov, \emph{Lectures on convex optimization}, second ed., Springer
  Optimization and Its Applications, vol. 137, Springer, Cham, 2018.
  \MR{3839649}

\bibitem{polyak1963gradient}
B.~T. Poljak, \emph{Gradient methods for minimizing functionals}, \v Z. Vy\v
  cisl. Mat i Mat. Fiz. \textbf{3} (1963), 643--653. \MR{158568}

\bibitem{rotaru2024exact}
Teodor Rotaru, Fran{\c{c}}ois Glineur, and Panagiotis Patrinos, \emph{Exact
  worst-case convergence rates of gradient descent: a complete analysis for all
  constant stepsizes over nonconvex and convex functions}, arXiv preprint
  arXiv:2406.17506 (2024).

\bibitem{taylor2017smooth}
Adrien~B. Taylor, Julien~M. Hendrickx, and Fran\c~cois Glineur, \emph{Smooth
  strongly convex interpolation and exact worst-case performance of first-order
  methods}, Math. Program. \textbf{161} (2017), no.~1-2, 307--345. \MR{3592780}

\bibitem{taylor2018exact}
\bysame, \emph{Exact worst-case convergence rates of the proximal gradient
  method for composite convex minimization}, J. Optim. Theory Appl.
  \textbf{178} (2018), no.~2, 455--476. \MR{3825634}

\bibitem{teboulle2023elementary}
Marc Teboulle and Yakov Vaisbourd, \emph{An elementary approach to tight worst
  case complexity analysis of gradient based methods}, Math. Program.
  \textbf{201} (2023), no.~1-2, 63--96. \MR{4620224}

\end{thebibliography}

\end{document}